\documentclass[12pt,onesided]{amsart}
\usepackage{caption}
\usepackage{subcaption}
\usepackage{graphicx}  %[dvips]
\usepackage[margin=0.6in,top=1.5in ]{geometry}
\usepackage{hyperref}
\hypersetup{
    colorlinks=true,
    linkcolor=blue,
    citecolor=red,
    urlcolor=blue,
    }
\graphicspath{{Images/}}

\usepackage{verbatim}

\usepackage{tikz}
\usetikzlibrary{arrows.meta}
\tikzset{%
  >={Latex[width=2mm,length=2mm]},
            base/.style = {rectangle, rounded corners, draw=black,
                           minimum width=4cm, minimum height=1cm,
                           text centered, font=\sffamily},
  activityStarts/.style = {base, fill=blue!30},
       startstop/.style = {base, fill=red!30},
    activityRuns/.style = {base, fill=green!30},
         process/.style = {base, minimum width=2.5cm, fill=orange!15,
                           font=\ttfamily},
}

\newtheorem{rmrk}{Remark}[section]

\newcommand{\R}{\mathbb{R}}

\newcommand{\N}{\mathbb{N}}

\newcommand{\xxi}{{\tilde \xi}}

\newcommand{\ma}{ {\mathcal{A}}}

\newcommand{\Id}{\textbf{\textup{Id}}}

\newcommand{\norm}{\lvert \lvert}

\newcommand{\cof}{\textup{cof}}

\newcommand{\tdet}{\textup{det}}

\newcommand{\md}{\mathcal{D}}

\newcommand{\ocg}{\overset{.} \gamma}

\begin{document}
%%-----------------------------
%%      the top matter
%%-----------------------------
\title[A First Insight into Solidification and Rupture of Blood]{Towards a Mathematical Model for the Solidification and Rupture of Blood in Stenosed Arteries}

%\thanks{...}\thanks{...}% At most 5 thanks
%
\author{Fatima Abbas}
\address{Laboratory of Mathematics (EDST)  - Lebanese University - Lebanon \\
Laboratoire de Math\'ematiques Appliqu\'ees du Havre (LMAH) - Le Havre University - France \\
Email: fatima-abs@hotmail.com}
\author{Ayman Mourad}
\address{Department of Mathematics - Faculty of Science (I) - Lebanese University - Lebanon \\
Email: ayman.mourad@ul.edu.lb}
%
%\date{January 20, 2020}
%
\begin{abstract} 
In this paper, we present a mathematical and numerical model for blood solidification and its rupture in stenosed arteries. The interaction between the blood flow and an existing stenosis in the arterial wall is modeled as a three dimensional fluid-structure interaction problem.
The blood is assumed to be a non-Newtonian incompressible fluid with a time-dependent viscosity that obeys a modified Carreau's model and the flow dynamics is described by the Navier-Stokes equations. Whereas, the arterial wall is considered a hyperelastic material whose displacement satisfies the quasi-static equilibrium equations. 
Numerical simulations are performed using FreeFem++ on a two dimensional domain. We investigate the behavior of the viscosity of blood, its speed and the maximum shear stress. From the numerical results, blood recirculation zones have been identified. Moreover, a zone of  blood of high viscosity and low speed has been observed directly after the stenosis in the flow direction. This zone may correspond to a blood accumulation and then solidification zone that is subjected to shear stress by the blood flow and to forces exerted by the artery wall deformation. Therefore, this zone is thought to break and then to release a blood clot that leads to the occlusion of small arterioles.

\end{abstract}
%
%\begin{resume} 
%Dans cet article, nous présentons un modèle mathématique et numérique pour la solidification du sang et sa rupture dans les artères en présence d'une sténose. L'interaction entre le flux sanguin et la st\'enose qui se trouve dans la paroie art\'erielle  est modélisée par un problème d'interaction fluide-structure en trois dimensions.
%Le sang est supposé être un fluide non-newtonien incompressible dont la viscosité dépend du temps et mod\'eli\'see par un modèle de Carreau modifi\'e. Son écoulement est régi par les équations de Navier-Stokes. La paroi artérielle est consid\'er\'ee comme un matériau hyperélastique dont le déplacement satisfait les équations d'\'equilibre quasi-statique.
%Des simulations numériques ont \'et\'e effectuées à l’aide du logiciel FreeFem++ sur un domaine bidimensionnel. Le comportement de la viscosité du sang et sa contrainte de cisaillement maximale ont \'et\'e \'etudi\'es. Les résultats num\'eriques ont montr\'e l'apparition d'une zone de recirculation du sang. De plus, une zone de sang ayant une viscosit\'e \'elev\'ee et une vitesse assez petite a \'et\'e observ\'ee. Cette zone pourrait correspondre \`a une zone d'accumulation et puis de solidification du sang et elle subit \`a des contraintes de cisaillement par le flux sanguin ainsi que \`a des forces exerc\'ees par la d\'eformation de la paroie art\'erielle. Par cons\'equent, ceci pourrait conduire \`a la rupture de cette zone et ensuite \`a la formation d'un caillot qui finit par l'occlusion de petites arterioles. 
% \end{resume}
%

\maketitle
%%-----------------------------
%%      your text
%%-----------------------------
\label{sec:into}

\section{Introduction}
Cardiovascular diseases, mainly due to atherosclerosis, are causes with the highest percentage leading to death worldwide.
Curiosity of finding cures for these diseases has enthused mathematicians to study them from their mathematical viewpoint. Consequently, many computational techniques and models have been developed. These tools aim at describing the blood flow and are effictive in studying the response of the arterial wall under certain conditions, the characteristics of the blood components 
%(Red Blood Cells, White Blood Cells and plasma) 
in addition to those of the heart \cite{Buriev,Peskin1,Peskin2,Tokarev1,Tokarev2,Yue}. 
Indeed, mathematical models with numerical analysis and simulations play an important role in providing knowledge and insights that are 
%indistinguishable 
unnoticeable clinically.
The rheological behavior of blood is captured by deriving constitutive models that constitute a constructive tool in the diagnosis of the pathologies, investigating appropriate remedies and proposing preventive therapies \cite{Barnes,Maco}. 
Recently, the lumen-wall modeling has been adopted using fluid-structure interaction (FSI) model. In the FSI, the behavior of the blood and the arterial wall are taken into consideration, so that one is capable of representing them by their appropriate dynamics and models. A most commonly used method when dealing with FSI systems is the Arbitrary Lagrangian-Eulerain (ALE) method \cite{Donea} that is effective when combining the fluid formulation in the Eulerian description and the structure formulation in the Lagrangian description.
An overview of FSI in biomedical applications has been considered in the book \cite{Tomas}. In particular, modeling of cardiovascular diseases has been highlighted in \cite{Luca,Quart,Zhan}.
In lumen-wall modeling some difficulties encounter
due to the complexity of the arterial wall formed of several layers, each with its own unique mechanics and thickness. 
Assuming that the arterial wall is negligibly thin, or the ratio of the arterial wall thickness to the aretry raduis is small, reduced shell or membrane models have been employed \cite{ani2006,Tallec}.
Introduction of computational model with FSI in order to investigate the wall shear stresses, blood flow field and recirculation zones in stenosed arteries have been studied in \cite{Buriev} where the blood is considered to be an incompressible Newtonian fluid, whereas in the case of a compressible non-Newtonian fluid these factors have been analyzed in \cite{Cho}. 
%For most FSI problems, analytical solutions are impossible to realize. 
%
%Nevertheless, approximate solutions can be obtained by employing numerical simulations. 
%
\\

Currently, numerous computational models are simple when describing the cardiovascular diseases and the related processes such as inflammations, coagulation, plaque growth, clot formation, etc.. Indeed, they are managed to capture only some essential features of the processes that take place in the cardiovascular system. Further, these models neglect some of the blood components with their characteristics as well as the arterial wall layers and their own mechanics. 
Consequently, more suitable models are needed through which the physiological parameters associated to the blood, atherosclerosis and clots must be investigated clinically.
In addition, the arterial wall must be considered as a multi-component structure taking into account the effects of the plaque growth and the clot formation on their mechanism. Further, timescales of the biological phenomena, the pulse duration and the time of the plaque growth must be analyzed.
\\

In the present work, the interaction between the blood modeled by a modified Carreau's model and the hyperelastic incompressible arterial wall has been considered.
In a first step, we introduce the FSI system which is composed of the incompressible Navier-Stokes equations representing the blood flow dynamics and the quasi-static equilibrium equations describing the elastic large deformation of the arterial wall. In addition, coupling conditions that ensure a global energy balance of the FSI system have been imposed on the common interface. Variational formulation has been presented and its discrete formulation has been derived. In fact, the Navier-Stokes equations have been semi-discretized in time, while, the nonlinear material equilibrium equations have been solved using the Newton-Raphson method.
Numerical simulations have been performed using FreeFem++ \cite{freefem}. A deep analysis has been made for better understanding of the behavior of the blood flow, the blood viscosity, the maximum shear stress and the recirculation zones. Based on the numerical results, a location where the blood is thought to accumulate and solidify has been identified. Finally, the factors affecting this zone have been investigated, in particular, this zone is subjected to forces exerted by the artery wall and the blood flow. Consequently, numerical simulations for the deformation of this zone have been performed in order to understand its rupture and thus the release of a clot that will lead to the occlusion of small arterioles.

\section{The Fluid-Structure Interaction Problem}
\label{FSI model}

The total domain $\Omega(t)$ representing the artery in the actual configuration at time $t>0$ is composed of two sub-domains $\Omega_f(t)$ and $\Omega_s(t)$ representing the lumen of the artery and the arterial wall, respectively.
\\

We denote by $\tilde{\Omega}_s$ the reference configuration of the structure of density $\tilde{\rho}_s$. Its deformation is described by the displacement field $ {\tilde \xi}_s : \tilde{\Omega}_s \times \mathbb{R}^+ \longrightarrow \mathbb{R}^3 $ that satisfies the quasi-static incompressible equilibrium equations. The evolution of the structure domain can also be given by the deformation map
$ \varphi_s : \tilde{\Omega}_s \times \mathbb{R}^+ \longrightarrow \mathbb{R}^3$ defined in terms of the displacement $ \tilde {\xi}_s$
as $ \varphi_s( \tilde {x},t)=  \tilde{x} +  \tilde {\xi}_s ( \tilde {x},t)$. 
Its deformation gradient $F_s :\tilde{\Omega}_s\times \mathbb{R}^+ \longrightarrow \mathbb{M}_3(\R)$ which is a second order tensor is given by $ F_s =  \nabla_{\tilde {x}} \varphi_s$. Its associated Jacobian is $J_s ( \tilde {x},t) = \tdet( F_s ( \tilde {x},t))$.
\\

On the other hand, we describe the blood flow dynamics by the incompressible Navier-Stokes equations on the sub-domain $\Omega_f(t)$. We denote by
\[
 v: \Omega_f(t) \times \mathbb{R}^+ \longrightarrow \mathbb{R}^3
\quad\textrm{and}\quad
p_f: \Omega_f(t) \times \mathbb{R}^+ \longrightarrow \mathbb{R}
\]
%%%
the velocity of the blood and its pressure, respectively. Further, the blood is assumed to be an incompressible fluid with a constant density $\rho_f$.
\\
The sub-domain $\Omega_f(t)$ of moving boundaries evolves from some reference configuration $\tilde{\Omega}_f$ according to an ALE map $\ma$ given by
\begin{align}\label{ALE-fluid-ch1}
\ma(.,t):\tilde{\Omega}_f &\longrightarrow \Omega_f(t) \nonumber \\
{ \tilde x} & \longrightarrow  \mathcal{A}({\tilde x},t)= x
\quad \textrm{for} \ t \in \mathbb{R}^+,
\end{align}
that is, $\Omega_f(t) =  \ma(\tilde{\Omega}_f,t)$.
\\
The ALE map $\mathcal{A}$ is considered to be an extension of the displacement $ {\tilde \xi}_s$ of the interface $\tilde{\Gamma}_c
=
\partial \tilde{\Omega}_s \cap \partial \tilde{\Omega}_f$, that is to say
\begin{align}\label{extension}
 \ma( {\tilde x},t)= {\tilde x} + \mathcal{E}xt( \tilde{\xi}_s( {\tilde x},t)|_{\tilde{\Gamma}_c}).
\end{align}
The operator $\mathcal{E}xt$ stands for an extension of the displacement of the boundary $\tilde{\Gamma}_c$. Possible extensions can be found in \cite[Section 5.3]{Richter}, \cite[Chapter 2]{Vincent} (harmonic, biharmonic, wislow, etc.). In particular we consider the harmonic extension as we will see in Subsection \ref{Discrete}. The deformation gradient associated to $ \ma$ is $ F_f :\tilde{\Omega}_f\times \mathbb{R}^+ \longrightarrow \mathbb{M}_3(\R)$ defined by $ F_f=  \nabla_{ {\tilde x}}  \ma $ where the symbol $ \nabla_{ {\tilde x}}$ indicates the gradient with respect to the variable $ {\tilde x}=(\tilde{x}_1,\tilde{x}_2,\tilde{x}_3)$.
 Its Jacobian is $J_f ( {\tilde x},t) = \tdet( F_f ( {\tilde x},t))$.
\\

Here and throughout the context, $ {\tilde \xi}_f$ denotes the displacement of the domain $\tilde{\Omega}_f$ which we set to be $\mathcal{E}xt( \widetilde{\xi}_s|_{\tilde{\Gamma}_c})$. Formulating the Navier-Stokes equations in the ALE frame results a new variable $ w$ that describes the velocity of the domain $\Omega_f(t)$. It is related to the displacement $ {\tilde \xi}_f$ by the relation $ w=\partial_t  {\tilde \xi}_f \circ  \ma^{-1}$. It is worth to point out that $ w \neq  v$. One must distinguish between $ v$ the physical velocity of the particles and $ w$ the velocity of the fluid domain $\Omega_f(t)$.

In what follows, we refer to the elements in the reference configuration by " $\widetilde{ }$ ". In fact the velocity and the pressure of the blood are given on the reference configuration $\tilde{\Omega}_f$ by
\begin{align}\label{elts in Ref. conf.}
 {\tilde v}( {\tilde x},t)
=
 v
\big(
 \ma( {\tilde x},t),t
\big)
\quad \textup{and} \quad
\tilde{p}_f ( {\tilde x},t)
=
p_f
\big(
 \ma( {\tilde x},t),t
\big)
\qquad \forall \ ( {\tilde x},t) \in \tilde{\Omega}_f \times \mathbb{R}^+.
\end{align}
The Cauchy stress tensor $ \sigma_f( v,p_f)$ is expressed in terms of the strain tensor
$ D( v)=\dfrac{ \nabla  v+ ( \nabla  v)^t }{2}$
as
\begin{align}\label{Cauchy-stress}
 \sigma_f( v,p_f)
=
2 \mu  D( v)- p_f \ \Id,
\end{align}
where $\mu=\mu(D(v))$ represents the blood viscosity that will be detailed in the sequel.
In the reference configuration $\tilde{\Omega}_f$, the stress tensor is given by
\begin{align*}
{\tilde{\sigma}}_f({\tilde{v}},\tilde{p}_f)
=
\mu
\big(
 \nabla {\tilde{v}} \ ( \nabla  \ma)^{-1}
+
( \nabla  \ma)^{-t} ( \nabla {\tilde{v}})^t
\big)
-\tilde{p}_f \ {\Id}.
\end{align*}
%%%
%%%
%%%

The arterial wall is assumed to be a hyperelastic material then it is characterized by the existence of an energy density function $W( F_s)$ such that the first Piola-Kirchhoff stress tensor $ P = \dfrac{\partial W( F_s)}{\partial  F_s}$. Further, due to the incompressible behavior of the material its Piola-Kirchhoff stress tensor is modified to the form  
\begin{align*}
{P}_\textup{inc}
%( {\tilde x})
={P}
%({\tilde x})
+ \tilde{p}_{hs} \cof( F_s).
\end{align*}
The variable $\tilde{p}_{hs}$, called the hydrostatic pressure, plays the role of the Lagrange multiplier associated to the incompressibility condition $\tdet( F_s)=1$.
\\

On the fluid domain $\Omega_f(t)$, a volumetric force $ f_f: \Omega_f \times \R^+ \longrightarrow \R^3 $ is applied. Moreover, a velocity $ v_\textup{in}$ is enforced on the inlet of the artery $\Gamma_\textup{in}(t)$. On the contrary, a free-exit condition given by $ \sigma_f( v,p_f)  n_f=0$ is enforced on the outlet $\Gamma_\textup{out}(t)$.
\\
On the other hand, a volumetric force $ f_s: \tilde{\Omega}_s \times \R^+ \longrightarrow \R^3 $ is applied on the structure domain which is assumed to be fixed on the boundary $\tilde{\Gamma}_2$, that is to say, $ {\tilde \xi}_s =0$ on $\tilde{\Gamma}_2$.
\\
On the interface $\Gamma_c(t)$, surface forces $ g_f: \Omega_f \times \R^+ \longrightarrow \R^3 $ and $ g_s: \Omega_s \times \R^+ \longrightarrow \R^3 $ are exerted from the fluid domain and the structure domain, respectively. \\

\begin{figure}[h!]
\centering
\includegraphics[scale=0.5]{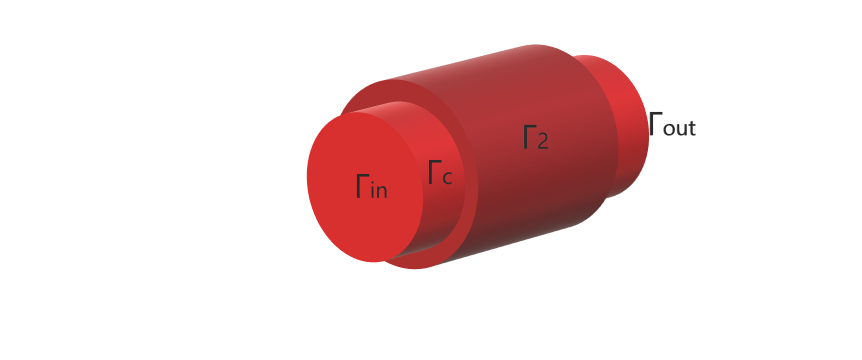}
\caption{A 3D model of an artery in the actual configuration.}
\label{3D model of artery}
\end{figure}

The FSI model describing the blood-wall interaction is obtained by the coupling 
between the incompressible Navier-Stokes equations which are formulated in the ALE frame and the quasi-static incompressible elasticity equations formulated in the Lagrangian frame on the reference configuration $\tilde{\Omega}_s$. The FSI system is\\
\\
Find 
\begin{center} 
${ \tilde{v} }: \tilde{\Omega}_f \times \R^+ \longrightarrow \R^3,\quad \tilde{p}_f : \tilde{\Omega}_f \times \R^+ \longrightarrow \R,
 \quad {   {\tilde \xi}_f }: \tilde{\Omega}_f \times \R^+ \longrightarrow \R^3$,
\end{center}
\begin{center}
$ {   {\tilde \xi}}_s : \tilde{\Omega}_s \times \R^+ \longrightarrow \R^3,\quad\tilde{p}_{hs} : \tilde{\Omega}_s \times \R^+ \longrightarrow \R$,
\end{center}
such that

\begin{equation}\label{FSI-Problem}
\begin{cases}
\rho_f \partial_t  v|_{\ma}
+ 
\rho_f ( v-  w)^t  \nabla_{ x}  v
-
 \nabla_{ x} \cdot  {\sigma}_f( v,p_f)
= \rho_f  {f}_f 
&
\qquad \textrm{on} \quad \Omega_f(t) \times (0,T),

\\

 \nabla_{ x} \cdot  v = 0
&
\qquad \textrm{on} \quad \Omega_f(t) \times (0,T),

\\
 v =  v_\textup{in} &
\qquad \textrm{on} \quad \Gamma_\textup{in}(t) \times (0,T),

\\

 \sigma_f( v,p_f)=0
&
\qquad \textrm{on} \quad \Gamma_\textup{out}(t) \times (0,T),

%\\
%
%
% \sigma_f( v,p_f)  n_f= g_f
%&
%\qquad \textrm{on} \quad \Gamma_c(t) \times (0,T),

\\
%~~~~~~~
%J_s \tilde{\rho}_s \partial^2_t \tilde{\xi}
-
 \nabla_{ {{\tilde x}}}  \cdot {P}_\textup{inc}( {{\tilde x}}) =J_s \tilde{\rho}_s {\tilde{f}}_s
&
\qquad \textrm{on} \quad \tilde{\Omega}_s \times (0,T),

%~~~~~~~
\\
J_s=1 & \qquad \textrm{on} \quad  \tilde{\Omega}_s \times (0,T),

%~~~~~~~
\\
 {\tilde{\xi}}_s = 0 & \qquad \textrm{on} \quad  \tilde{\Gamma}_{2}\times (0,T),

%~~~~~~~
%\\
%{P}_\textup{inc}({{\tilde x}}) {\tilde{n}}_s = J_s  {\tilde{g}}_s 
%
%& \qquad \textrm{on} \quad  \tilde{\Gamma}_c\times (0,T),

%~~~~~~~~
\\

 v = w 
& \qquad \textrm{on} \quad  \Gamma_c(t)  \times (0,T),

%~~~~~~~~~

\\

{P}_\textup{inc}( {\tilde x}) 
{\tilde{n}}_s + J_f {\tilde{\sigma}}_f({\tilde{v}},\tilde{p}_f)  F^{-t}_f {\tilde{n}}_f
=0
& \qquad \textrm{on} \quad  \tilde{\Gamma}_c \times (0,T) ,
\end{cases}
\end{equation}
where $ {\tilde v}$ and $\tilde{ p}_f$ are given by \eqref{elts in Ref. conf.}. Further, $\tilde{\Gamma}_c$ is the transformation of $\Gamma_c$ to the reference configuration. Figure \ref{3D model of artery} illustrates a 3D model of an artery in the actual configuration including the boundaries.
\begin{rmrk}\label{rel btwn extension and deformation}
From Expression \eqref{extension} we get that the ALE map $\ma$ and the structure deformation $ \varphi_s$ coincide on the interface $\tilde{\Gamma}_c$, that is to say, 
\begin{align*}
 \varphi_s \equiv  \ma \qquad \textup{on} \ \tilde{\Gamma}_c.
\end{align*}
\end{rmrk}
\begin{rmrk}
Due to the incompressibility condition the $ij$-th component of $ \sigma_f( v,p_f)$ is 
\[
 \sigma_{ij}=-p_f \delta_{ij} 
 +
 \mu \bigg(
 \dfrac{\partial v_i}{\partial x_j}+\dfrac{\partial v_j}{\partial x_i}
 \bigg), \qquad i,j=1,2,3,
\]
where $ \delta_{ij}$ is the Kronecker symbol. 
The shear stress components are $\sigma_{12},\ \sigma_{13}$ and $ \sigma_{23}$, whereas $\sigma_{11},\ \sigma_{22}$ and $\sigma_{33}$ are the normal stress components.

In a two dimensional space the maximum shear stress- an effective parameter in studying the forces exerted on a fluid- is given by the expression \textup{\cite{Yue}}
\begin{equation}\label{max-shs}
\sigma_{max} = \sqrt{\left(\dfrac{\sigma_{11}-\sigma_{22}}{2}\right)^2 + \sigma_{12}^2}.
\end{equation}
\end{rmrk}
The variational formulation associated to System \eqref{FSI-Problem} is 
\\
Find 
\begin{center} 
${ \tilde{v} }: \tilde{\Omega}_f \times \R^+ \longrightarrow \R^3$, 
$\tilde{p}_f : \tilde{\Omega}_f \times \R^+ \longrightarrow \R$,
 $ {   {\tilde \xi}_f }: \tilde{\Omega}_f \times \R^+ \longrightarrow \R^3$,
\end{center}
\begin{center}
$ {   {\tilde \xi}}_s : \tilde{\Omega}_s \times \R^+ \longrightarrow \R^3$,
$ \tilde{p}_{hs} : \tilde{\Omega}_s \times \R^+ \longrightarrow \R$,
\end{center}
such that
\begin{equation}\label{vf1}
\begin{aligned}
{\tilde \xi}_f &= \mathcal{E}xt( {\tilde \xi}_s|_{\tilde{\Gamma}_c}) \ \text{and} \ 
 {\tilde w}= \dfrac{\partial  \tilde{\xi}_f }{\partial t}
&& \text{in} \quad \tilde{\Omega}_f, \\
 v & =  w && \text{on} \quad \Gamma_c(t), \\
 {\tilde \xi}_s & = 0 && \text{on} \quad \tilde{\Gamma}_2.  
\end{aligned}
\end{equation}
and
\begin{equation}\label{vf2}
\begin{cases}
\rho_f
\displaystyle\int_{\Omega_f(t)} \dfrac{\partial  v}{\partial t}\Big|_{\ma} \cdot  \eta_f \ \ d  x
+
\rho_f \int_{\Omega_f(t)}  ( v- w)^t  \nabla_{ x}  v \cdot  \eta_f \ d  x
+
\int_{\Omega_f(t)} {\sigma}_f( v,p_f):  \nabla_{ x}  \eta_f \ d  x
\\
-
\displaystyle\int_{\Gamma_c(t)}
 {\sigma}_f( v,p_f) n_f \cdot  \eta_f \ d\Gamma
%%%
%%%
%%
=
\rho_f \int_{\Omega_f(t)}   f_f \cdot  \eta_f \ d  x,
\\
\displaystyle\int_{\Omega_f(t)} q_f \ \nabla_{ x} \cdot  v \ d  x =0,
\end{cases}
\end{equation}
\begin{equation}\label{Variational form-Elastodynamic}
\begin{cases}
%\rho_s
%\displaystyle\int_{\tilde{\Omega}_s}
%\partial_t^2   {\tilde \xi}  {\tilde \eta}_s \ d\tilde{x}
%+
%%
%%
\displaystyle \int_{\tilde{\Omega}_s} {P}
%({\tilde x}) 
:  \nabla_{ {\tilde x}}  {\tilde \eta}_s \ d  {\tilde x}
+
\int_{\tilde{\Omega}_s} \tilde{p}_{hs} \cof( F_s)
:  \nabla_{ {\tilde x}}  {\tilde \eta}_s \ d {\tilde x}
\\
-
\displaystyle\int_{\tilde{\Gamma}_c}  {P}_\textup{inc}( {\tilde x}) 
{\tilde{n}}_s \cdot   {\tilde \eta}_s \
d{\tilde \Gamma}

=
\int_{\tilde{\Omega}_s} J_s \rho_s  {\tilde f}_s \cdot  {\tilde \eta}_s \ d  {\tilde x}. 

\\

\displaystyle\int_{\tilde{\Omega}_s} \tilde{q}_s  \ (J_s-1) \ d {\tilde x}=0 
\end{cases}
\end{equation}
for all 
 $( \eta_f,q_f) \in H^1_{\Gamma_{in}(t)}(\Omega_f(t)) \times L^2(\Omega_f(t))$  
 and
 $({\tilde \eta}_s,\tilde{q}_s) \in H^1_{\tilde{\Gamma}_2}(\tilde{\Omega}_s) \times L^2(\tilde{\Omega}_s)$.
%where $ v = {\tilde v} \circ  \ma^{-1}_t$ and 
%$p_f = \tilde{p}_f \circ  \ma_t^{-1}$.
%%
\\
The coupling conditions on the interface $\tilde{\Gamma}_c$ are given in the strong form as

\begin{equation}\label{vf4}
\begin{cases} 
 v \circ  \ma = \partial_t  {\tilde \xi} , \\
{P}_\textup{inc}( {\tilde x}) 
{\tilde{n}}_s + J_f {\tilde{\sigma}}_f({\tilde{v}},\tilde{p}_f)  F^{-t}_f {\tilde{n}}_f
=0.
\end{cases}
\end{equation}

The spaces $H^1_{\Gamma_{in}(t)}(\Omega_f(t))$ and $H^1_{\tilde{\Gamma}_2}(\tilde{\Omega}_s)$ are respectively
\begin{align*}
H^1_{\Gamma_{in}(t)}(\Omega_f(t))=
\left\{
 \eta_f \in H^1(\Omega_f(t))
; \ \eta_f=0 \ \ \textup{on} \ \ \Gamma_\textup{in}(t) 
\right\}
\
\textrm{and}
\quad
H^1_{\tilde{\Gamma}_2}(\tilde{\Omega}_s)
=
\Big\{
\tilde {\eta}_s \in H^1(\tilde{\Omega}_s); \ \tilde {\eta}_s = 0 \ \ \textup{on} \ \ \tilde{\Gamma}_2 \Big\}.
\end{align*}
\subsection{The Discrete Variational Formulation of the FSI Problem}\label{Discrete}
The variational formulation \eqref{vf1}-\eqref{vf4} of the FSI stands for the incompressible homogeneous Navier-Stokes equations coupled with the quasi-static incompressible equilibrium equations. We assume that no external forces are exerted on neither the fluid domain nor the structure domain, i.e, $ f_f=0$ and $ {\tilde f}_s=0$. Consider a time step $\Delta t >0$ and finite element partitions $\mathcal{V}_h$ and $\mathcal{W}_h$ for the fluid and the solid sub-domains respectively, of a maximum diameter denoted by $h$.  
Our aim is to approximate the solution $( v,p_f, {\tilde \xi}_s, \tilde{p}_{hs},\ma)$ at time $t_n=n \Delta t$, for $n \in \N$, in the finite element spaces. The approximation of the solution at time $t_n$ is denoted by $( v^n,p_f^n, {\tilde \xi}_s^n,\tilde{p}^n_{hs} ,\ma^n)$. 
%%%%
%%%%
%%%%
%%%%
\subsubsection*{Semi-Discretization in Time of the Fluid Sub-Problem}
In order to guarantee the existence and uniqueness of the solution of the discrete fluid sub-problem when performing the numerical simulations, we use the penalty method \cite{Glo}.
This method consists of replacing the natural weak formulation by a regular one by adding a term multiplied by a sufficiently small parameter $\epsilon \ll 1$. 
Indeed, writing the modified formulation in a matrix form results a positive definite matrix, which assures the existence and the uniqueness of the solution of the discrete sub-problem. 
%Moreover, the solution of the new formulation converges to the solution of the original problem.
%
The weak formulation associated to the Navier-Stokes equations obtained upon adding a negligible parameter $\epsilon$ is then semi-discretized in time, that is, the convective term and the viscosity are considered at the instant $t_n$, whereas other terms are considered at time $t_{n+1}$.
The discrete formulation reads
\begin{equation}\label{discrete-v.f-fluid}
\begin{split}
&
\rho_f
\dfrac{1}{\Delta t} 
\displaystyle
\int_{\Omega_f(t_n)}  v^{n+1} \cdot  \eta_f \ d x
+
\rho_f
\dfrac{1}{\Delta t}
\int_{\Omega_f(t_n)}
( v^n \circ  X^n)
\cdot \  \eta_f \ \ d  x
-
\rho_f 
\int_{\Omega_f(t_n)}
( w^{n+1})^t   \nabla_{ x}  v^{n+1} \cdot  \eta_f \ d  x
\\
&
+
2\int_{\Omega_f(t_n)} \mu^{n}  D(  v^{n+1})  :  \nabla_{ x}  \eta_f \ d x
-
\int_{\Omega_f(t_n)} p_f^{n+1} \ \nabla_{ x} \cdot  \eta_f \ d  x
-
2 
\int_{\Gamma_c(t_n)}
\mu^n D(v^{n+1}) n_f \cdot  \eta_f \ d\Gamma 
\\
&
+
\int_{\Gamma_c(t_n)}
 p_f^{n+1} n_f  \cdot  \eta_f \ d\Gamma 
+
\int_{\Omega_f(t_n)} q_f  \ \nabla_{ x} \cdot  v^{n+1} \ d  x
+
\int_{\Omega_f(t_n)} \epsilon \ p_f^{n+1}q_f \ d  x
=
0, 
\end{split}
\end{equation}
where the non-linear convective term $
\dfrac{1}{\Delta t}
( v^n \circ  X^n)
$ can be approximated by  \cite[Section 9.5, p. 267]{freefem}
\begin{align}\label{approx}
\dfrac{1}{\Delta t}
\Big[
 v \big( x- v(x,t_n)\Delta t,t_n \big)
\Big] .
\end{align}

Notice that, since the strain rate tensor $ D(v)$ is symmetric, then we have $ D( v): \nabla  \eta_f =  D( v): D( \eta_f)$ which gives 
\begin{align*}
\int_{\Omega_f(t_n)} \mu^{n}  D(  v^{n+1})  :  \nabla_{ x}  \eta_f \ d  x
=
\int_{\Omega_f(t_n)} \mu^{n}  D(  v^{n+1})  :  D ( \eta_f )\ d  x.
\end{align*}

\subsubsection*{Newton-Raphson Method for the Structure Sub-Problem}

Regarding the structure sub-problem, at the time iteration $t_{n+1}$ we will solve the non-linear problem \eqref{Variational form-Elastodynamic} using Newton-Raphson method. The variational formulation corresponding to the structure sub-problem at the iteration $t_{n+1}$ is 
\begin{equation}
\label{Variational form-Elastodynamic- iteration n+1}
\begin{cases}
\displaystyle
  \int_{\tilde{\Omega}_s} 
{P}^{n+1} 
: 
 {\nabla \tilde{\eta}_s}
 \ d {\tilde x}

+
 \int_{\tilde{\Omega}_s} 
\tilde{p}^{n+1}_{hs} 
{\cof(F}^{n+1}_s)
: 
{ \nabla \tilde{\eta}_s} \ d  {\tilde x}

\\
\displaystyle
%-
% \int_{\tilde{\Gamma}_c} 
%J_s^{n+1} 
% {\tilde g}^{n+1}_s 
%\cdot 
%{ \tilde{\eta}_s}  \ d\tilde{\Gamma} 

-
 \int_{\tilde{\Gamma}_c} 
P^{n+1} 
 {\tilde n}_s 
\cdot 
{ \tilde{\eta}_s}  \ d\tilde{\Gamma} 

-
\int_{\tilde{\Gamma}_c} 
\tilde{p}^{n+1}_{hs} 
{\cof(F}^{n+1}_s)
 {\tilde n}_s 
\cdot 
{ \tilde{\eta}_s}
\ d\tilde{\Gamma}

-
 \int_{\tilde{\Omega}_s}  
J^{n+1}_s 
 \tilde{\rho}_s 
  {\tilde f}_s 
 \cdot 
 { \tilde{\eta}_s}
  \ d  {\tilde x}

=0

\quad  \forall \ { \tilde{\eta}_s} \in \tilde{V}_s,

\\
\displaystyle

 \int_{\tilde{\Omega}_s} \tilde{q}_s (J_s^{n+1}-1) d {\tilde x} = 0,

\qquad \forall \ \tilde{q}_s \in L^2(\tilde{\Omega}_s).

\end{cases}
\end{equation}

The method depends on linearizing the structure sub-problem \eqref{Variational form-Elastodynamic- iteration n+1} with respect to the unknowns $ \varphi_s$ and $\tilde{p}_{hs}$.
We start initialization by considering a suitable choice of the initial values ($ \varphi_{s,0},\tilde{p}_{hs,0}$). In particular, we link the iterations of the Newton-Raphson method with the time iteration $t_n$ by considering $ \varphi_{s,0}= \varphi_s^n$ and $\tilde{p}_{hs,0}=\tilde{p}_{hs}^n$. Then, we solve iteratively the obtained system corresponding to the Newton-Raphson method until its solution converges to a solution of the non-linear System \eqref{Variational form-Elastodynamic}. To ensure the existence of the solution of the structure problem \eqref{Variational form-Elastodynamic- iteration n+1} we use the penalty method by modifying System \eqref{Variational form-Elastodynamic- iteration n+1} through adding 
the penalized term $\epsilon \int\limits_{\tilde{\Omega}_s} \tilde{p}_{hs} \tilde{q}_s\ d {\tilde x}$ with $\epsilon \ll 1$. For simplicity of notation, in what follows we omit the subscript $s$ of the deformation $ \varphi_s$, that is, we write $ \varphi_s \equiv  \varphi$. 
\\

We proceed to derive the formulation of the structure sub-problem corresponding to the Newton-Raphson method. Let us define the following space
\[
\mathcal{Z}=\lbrace
 \varphi=(\varphi_1,\varphi_2,\varphi_3):\tilde{\Omega}_s \longmapsto \R^3
,  \varphi=  \varphi^n \ \textup{on} \ \tilde{\Gamma}_c \ \ \textup{and} \;  \varphi=0 \; \textup{on} \; \tilde{\Gamma}_2 
\rbrace.
\]
Given $N \in \N$, a tolerance $tol$ and
\[
( \varphi_0,\tilde{p}_{hs,0}) \in \mathcal{Z} \times L^2(\tilde{\Omega}_s),
\]
we construct iteratively the two sequences $( \varphi_k)_{k \geq 1}$ and $(\tilde{p}_{hs,k})_{k \geq 1}$ by solving for $(\delta  \varphi_k,\delta \tilde{p}_{hs,k})$ the following system:
\\
\\
Set  $ \varphi_0= \varphi^n$.
Repeat: for $0 \leq k \leq N$, while $\norm \delta  \varphi_k \norm_2 \geq tol$, find ($\delta  \varphi_k$,$\delta \tilde{p}_{hs,k}$) in $\mathcal{Z} \times L^2(\tilde{\Omega}_s)$ satisfying
\\
\begin{equation}\label{linearizedVf}
\begin{cases}

\displaystyle 
 \int_{\tilde{\Omega}_s} 
\dfrac{\partial  P}{\partial  {F}_s}
( \nabla_{{\tilde x}}  \varphi_k) 
 \nabla_{{\tilde x}} \delta  \varphi_k :  \nabla_{{\tilde x}} {\tilde \eta}_s \ d  {\tilde x}

+
{\int_{\tilde{\Omega}_s}} 
\tilde{p}_{hs,k} \dfrac{\partial \cof}{\partial  F_s}( \nabla_{{\tilde x}}  \varphi_k) 
 \nabla_{{\tilde x}} 
\delta  \varphi_k :
 \nabla_{{\tilde x}} { \tilde \eta}_s 
\ d  {\tilde x}

\\ \\

\displaystyle
+{\int_{\tilde{\Omega}_s}} 
\delta \tilde{p}_{hs,k} 
\ \cof( \nabla_{{\tilde x}}  \varphi_k) :
 \nabla_{{\tilde x}}  {\tilde \eta}_s 
\ d{\tilde x}

-
{\int_{\tilde{\Gamma}_c}} 
 {\tilde \sigma}_{s,k} ( {\tilde x}) 
\dfrac{\partial \cof}{\partial  F_s}( \nabla_{{\tilde x}}  \varphi_k) 
: 
 \nabla_{{\tilde x}} \delta  \varphi_k \  {\tilde n}_s \cdot 
 {\tilde \eta}_s 
\ d\tilde{\Gamma}

\\ \\

\displaystyle
+ \int_{\tilde{\Omega}_s} 
{P}( \nabla_{{ \tilde x}}  \varphi_k):  \nabla_{{\tilde x}}  {\tilde \eta}_s \ d {\tilde x}
+
 \int_{\tilde{\Omega}_s} 
\tilde{p}_{hs,k} \ 
\cof ( \nabla_{{ \tilde x}}  \varphi_k):  \nabla_{{\tilde x}} {\tilde \eta}_s \ d {\tilde x}

\\ \\ 
\displaystyle

- \int_{\tilde{\Gamma}_c}
  {\tilde \sigma}_{s,k} ( {\tilde x}) 
 \  {\tilde n}_s \cdot 
  {\tilde \eta}_s \ d\tilde{\Gamma}

+
\epsilon
\int_{\tilde{\Omega}_s}
\tilde{p}_{hs,k} \ \tilde{q}_s \ d\tilde{x}

=
0 

\\ \\

\displaystyle
{\int_{\tilde{\Omega}_s}} 
\tilde{q}_s \ 
\cof( \nabla_{{\tilde x}}  \varphi_k) 
:
 \nabla_{{\tilde x}} 
\delta  \varphi_k \ d {\tilde x}
+
{\int_{\tilde{\Omega}_s}} 
\tilde{q}_s (\tdet( \nabla_{{\tilde x}}  \varphi_k)-1)\ d {\tilde x}
=
0
\end{cases}
\end{equation}
for all $({\tilde \eta}_s,\tilde{q}_s) \in \tilde{V}_s \times L^2(\tilde{\Omega}_s)$.

Set $ \varphi_{k+1}= \delta  \varphi_k +  \varphi_k$ and $k=k+1$.
\\

When the condition $\norm \delta  \varphi_k \norm_2 < tol$ is fulfilled then convergence of the Newton-Raphson method is achieved. Thus, the solution of the structure sub-problem \eqref{Variational form-Elastodynamic- iteration n+1} is given by $ \varphi^{n+1}= \varphi_k$, for the last value of $k$ for which the Newton-Raphson method converges.
\subsubsection*{Space Discretization of the FSI Problem}
Space discretization of the variational formulation is carried out using the finite element method (FEM) \cite{Girault-Raviart}. We consider the two finite element spaces associated to the fluid weak formulation
\begin{align*}
V^f_h \subset V_f \ \textrm{and} \ W^f_h \subset L^2(\Omega_f)
\end{align*}
and those associated to the structure weak formulation
\begin{align*}
\tilde{V}^s_h \subset \tilde{V}_s
 \ \textrm{and} \ \tilde{W}^s_h \subset L^2(\tilde{\Omega}_s)
\end{align*}
where $V^f_h$, $W^f_h$, $\tilde{V}^s_h$ and $\tilde{W}_h^s$ are finite dimensional subspaces.
The functional spaces associated to the velocity and displacement fields are considered to be $P_2$, whereas those associated to the pressures (fluid and hydrostatic) are considered to be $P_1$.
In what follows, all terms are discretized in space as mentioned above, so that the approximation of solution in finite element spaces is $( v_h,p^f_h, {\tilde \xi}^s_h,\tilde{p}^{hs}_h, \ma_h)$ verifying \eqref{discrete-v.f-fluid} and \eqref{linearizedVf}.
%, though, for simplicity the subscript $h$ is omitted in the context.
\\
\\
Finally, the discrete variational formulation reads:
\vspace{2mm}

Given $( v^n_h,p_{f,h}^n, {\tilde \xi}_{s,h}^n,\tilde{p}_{hs,h}^n, \ma^n_h)$ and a tolerance $tol$,
 find $( v^{n+1}_h,p_{f,h}^{n+1}, {\tilde \xi}_{s,h}^{n+1},\tilde{p}_{hs,h}^{n+1}, \ma^{n+1}_h)$ such that
\begin{equation}\label{FSI-eq1}
\begin{cases}
 \ma^{n+1}_h=  {\tilde x}+ \mathcal{E}xt( {\tilde \xi}_{s,h}^{n+1}|_{\tilde{\Gamma}_c}) &
\textrm{in} \ \tilde{\Omega}_f
\\
 {\tilde w}^{n+1}_h

=
\dfrac{1}{\Delta t}( {\tilde \xi}_{f,h}^{n+1}- {\tilde \xi}_{f,h}^n) \simeq \partial_t  {\tilde \xi}_{f,h}^{n+1}
& \textrm{in} \ \tilde{\Omega}_f
,
 \\
%
%
%v^{n+1} 
%= 
 w^{n+1}_h
=
\partial_t  \ma^{n+1}_h \circ ( \ma^{n+1}_h)^{-1}
%w(\tilde{\eta}_f^{n+1})
& \text{on} \  \ma_h^{n+1}(\tilde{\Gamma}_c), \\
 {\tilde \xi}_{s,h}^{n+1} = 0 & \textup{on} \ \tilde{\Gamma}_2
\end{cases}
\end{equation}
\begin{equation}\label{FSI-eqn2}
\begin{split}
&\rho_f 
\dfrac{1}{\Delta t}
\displaystyle
\int_{\Omega_f(t_n)}
 v_h^{n+1} 
\cdot \  \eta^f_h \ \ d  x
+
\rho_f
\dfrac{1}{\Delta t}
\int_{\Omega_f(t_n)}
( v^n_h \circ  X^n_h)
\cdot \  \eta^f_h \ \ d  x
%
%
%+
%
%
-
\rho_f 
\int_{\Omega_f(t_n)}   (w^{n+1}_h)^t   \nabla_{ x}  v_h^{n+1} \cdot 
 \eta^f_h \ d  x
%
%
%\nonumber 
\\
&+
2\int_{\Omega_f(t_n)} \mu_h^{n}  D(  v_h^{n+1})  :  \nabla_{ x}  \eta_h^f \ d x
-
\int_{\Omega_f(t_n)} p_{f,h}^{n+1} \ \nabla_{ x} \cdot  \eta_h^f \ d  x
-
\int_{\Gamma_c(t_n)}
 g_{f,h}^{n+1} \cdot  \eta^f_h \ d\Gamma
%\nonumber 
\\
&+
\int_{\Omega_f(t_n)} 
q_h^f \ 
\nabla_{ x} \cdot  v_h^{n+1} \ d  x
%%%
%%%
%%%
+
\epsilon \int_{\Omega_f(t_n)} p^{n+1}_{f,h} q^f_h \ d  x
=0 \qquad \qquad \forall \ ( \eta^f_h,q^f_h) \in V_h^f \times W_h^f.
\end{split}
\end{equation}
where the non-linear convective term 
$
\dfrac{1}{\Delta t}
( v^n_h \circ  X^n_h )
$ is approximated by the Expression \eqref{approx}.
\\
The coupling conditions on $\tilde{\Gamma}_c$ are
\begin{equation}\label{discretized coupling cdtn}
\begin{cases}

 {\tilde \sigma}^{n+1}_{s,h} \ \tilde{n}_s
=
 -\bigg( 
\big( 
 2
\mu^n D(v_h^{n+1})-p_{f,h}^{n+1} \ \Id 
\big)
 n_f \bigg)\circ  \varphi^n_h
,  \\

\partial_t  {\tilde \xi}_{s,h}^{n+1}
=
 v^{n+1}_h \circ  \varphi^n_h.

\end{cases}
\end{equation}

Fix $N \in \N$. 
Set $ \varphi_{0,h}= \varphi^{n}_h$ and $ {\tilde \sigma}_{s,k}^h
=  {\tilde \sigma}_{s,h}^{n+1}$. Repeat: for $0 \leq k \leq N$, find ($\delta  \varphi_{k,h}$,$\delta \tilde{p}^h_{hs,k}$) in $\mathcal{Z} \times L^2(\tilde{\Omega}_s)$ satisfying
\begin{equation}\label{FSI-eqn3}
\begin{cases}
\displaystyle
 \int_{\tilde{\Omega}_s} 
\dfrac{\partial {P}}{\partial  F_s}(\nabla_{ {\tilde x}}  \varphi_{k,h}) 
 \nabla_{ {\tilde x}} \delta  \varphi_{k,h} :  \nabla_{ {\tilde x}}  {\tilde \eta}^s_h 
\ d{\tilde x}
+
 \int_{\tilde{\Omega}_s} 
\tilde{p}_{hs,k} ^h
\dfrac{\partial \cof}{\partial  F_{s,h}}
( \nabla_{ {\tilde x}}  \varphi_{k,h}) 
 \nabla_{ {\tilde x}} \delta  \varphi_{k,h} :  \nabla_{{\tilde x}}  {\tilde \eta}^s_h
\ d {{\tilde x}}
\\ \\
\displaystyle
+ 
\int_{\tilde{\Omega}_s} 
\delta \tilde{p}_{hs,k}^h \ \cof( \nabla_{{\tilde x}}  \varphi_{k,h}) 
: 
 \nabla_{ {\tilde x}} {\tilde \eta}^s_h \ d {\tilde x}
-
\int_{\tilde{\Gamma}_c} 
 {\tilde \sigma}^h_{s,k}({{\tilde x}}) 
\dfrac{\partial \cof}{\partial  F_{s,h}}
( \nabla_{{\tilde x}}  \varphi_{k,h}) 
:  \nabla_{{\tilde x}} 
\delta  \varphi_{k,h}  
\  {\tilde n}_s 
\cdot 
 {\tilde \eta}^s_h
\ d\tilde{\Gamma}
\\ \\
\displaystyle
+
\int_{\tilde{\Omega}_s}  
{P}( \nabla_{{\tilde x}}  \varphi_{k,h}):  \nabla_{{\tilde x}}  {\tilde \eta}^s_h
 \ d {{\tilde x}}
+

\int_{\tilde{\Omega}_s}  
\tilde{p}_{hs,k} 
\ 
\cof ( \nabla_{{\tilde x}}  \varphi_{k,h})
:  \nabla_{ {\tilde x}} 
 {\tilde \eta}^s_h\ d {{\tilde x}}
\\ \\

\displaystyle
-
 \int_{\tilde{\Gamma}_c} 
 {\tilde \sigma}^h_{s,k} ({{\tilde x}}) 
\ \cof ( \nabla_{{\tilde x}}  \varphi_{k,h}) 
 {\tilde n}_s 
\cdot 
 {\tilde \eta}^s_h
\ d\tilde{\Gamma}
+
\epsilon
 
\int_{\tilde{\Omega}_s} \tilde{p}^h_{hs,k} \ \tilde{q}^s_h \ d{{\tilde x}}
=
0 \qquad \qquad \forall \ ({\tilde \eta}^s_h,\tilde q^s_h) \in \tilde{V}_h^s \times \tilde{W}_h^s,

\\ \\

\displaystyle

\int_{\tilde{\Omega}_s}  \tilde{q}^s_h \ \cof( \nabla_{{\tilde x}}  \varphi_{k,h}) :  \nabla_{{\tilde x}} \delta  \varphi_{k,h} \ d{{\tilde x}}
+

\int_{\tilde{\Omega}_s}  \tilde{q}^s_h (\tdet( \nabla_{{\tilde x}}  \varphi_{k,h})-1)\ d{{\tilde x}}
=
0 \qquad \qquad \forall \ \tilde q^s_h \in \tilde{W}_h^s,
\end{cases}
\end{equation}
as long as the error $\norm \delta  \varphi_{k,h} \norm_2 \geq tol$, set $\varphi^{k+1}= \delta  \varphi_{k,h} +  \varphi_{k,h}$ and $k=k+1$.
If $\norm \delta  \varphi_{k,h} \norm_2 < tol$, then convergence of the Newton-Raphson method is achieved. Thus the deformation of the structure domain is given by $ \varphi^{n+1}_h= \varphi_{k,h}$, for the last value of $k$ for which the convergence is achieved. 
Whence at the time iteration $t_{n+1}$ the displacement of the structure domain is 
$ {\tilde \xi}_{s,h}^{n+1}= \varphi^{n+1}_h -  {\tilde x}$. 
\\

At the time iteration $t=t_{n+1}$, the ALE map is given by the following relation
\begin{align}\label{harmonic extension}
 \ma_h^{n+1} ({\tilde x},t)= {\tilde x} +  {\tilde \xi}_{f,h}^{n+1}({\tilde x},t).
\end{align}
where 
$
{\tilde \xi}_{f,h}({\tilde x},t)$ is constructed
using the harmonic extension \cite[Section 5.3, pp. 247]{Richter}, \cite[Chapter 2, pp. 54]{Vincent} of the displacement  $ {\tilde \xi}_{s,h}({\tilde x},t)$ of the boundary $\tilde{\Gamma}_c$.

As a result,
The fluid domain  and the structure domain evolve from their reference configuration according to
\begin{align*}
\Omega_f^{n+1}=  \ma^{n+1}(\tilde{\Omega}_f)
\qquad 
\text{and}
\qquad
\Omega^{n+1}_s =  \varphi^{n+1} (\tilde{\Omega}_s),
\end{align*}
respectively.
\subsection{The Algorithm}
The FSI problem is solved using the finite element software FreeFem++. The steps of the algorithm used to solve the FSI problem are stated below. Complexity of the algorithm is dependent on the triangulation of the mesh.
Literally,
\begin{enumerate}
\item[1-] At the time step $t_{n+1}$, we solve the Navier-Stokes equations on the domain ${\Omega}_f(t_n)$ in the ALE frame to find the velocity of the fluid $ v^{n+1}$ and its pressure $p^{n+1}_f$.
\item[2-] Using the continuity of stresses \eqref{discretized coupling cdtn}, we are able to get the boundary condition on $\tilde{\Gamma}_c$ expressed in terms of 
$ \sigma_f( v^{n+1},p^{n+1}_f)$.
\item[3-] Solve the quasi-static incompressible elasticity equations on the reference configuration $\tilde{\Omega}_s$ using Newton-Raphson method. Check the convergence test of the Newton-Raphson method. When convergence is achieved, the deformation $ \varphi^{n+1}$ is set to be the  solution of the last Newton's iteration $k$, consequently the deformation $ {\tilde \xi}_s^{n+1}$ is obtained. Thus we can proceed to get the ALE map $ \ma^{n+1}$ using the harmonic extension given by Equation \eqref{harmonic extension}.
\item[4-] Move the fluid domain using the map $ \ma^{n+1}$, and the structure domain using its deformation $ \varphi^{n+1}$ and proceed to the iteration $t_{n+2}$, then start again from step (1-), and so on.
\end{enumerate}
\section{Numerical Results}\label{Numerical Results}
In this section we present the numerical results concerning the blood flow through stenosed arteries after performing simulations over a defined interval of time. The study is done by solving System \eqref{FSI-eq1}-\eqref{FSI-eqn3} on a two dimensional domain representing the artery using the software FreeFem++.

Our work is concerned in analyzing variables including the speed, the viscosity and the wall shear stress of blood in a stenosed artery. Further, we intent to locate the recirculation zones in the lumen.
%Investigating the manner of some variables during the first second would give us an idea of what would happen on the predefined interval.
~In our work we consider the following numerical values
\[
 \rho_f=1.056 \ \text{g/cm}^3 \ \text{and} \ \Delta t = 10^{-2} \ \text{s}.  
\]
The blood is considered to be of a non-Newtonian behavior.
Its viscosity is assumed to obey Carreau model \cite{Adelia}
\begin{align}\label{Carreau}
\mu(\ocg)=
\mu_\infty +(\mu_0 - \mu_\infty)
\big[1+ (\lambda \ocg)^2 \big]^{\frac{n-1}{2}},
\end{align}
where $\ocg$ stands for the shear rate defined as
\begin{align}\label{shear rate}
\overset{.}\gamma
=
\sqrt{2 \ \textrm{tr}( D( v))^2}= \sqrt{-4 I_2}.
\end{align}
The parameters present in the Carreau model are
\[ \lambda = 3.313 \ \text{s},\quad n=0.3568,\quad \mu_\infty= 0.00345 \ \text{Pa.s} \quad \text{and} \quad \mu_0 = 0.056 \  \text{Pa.s}.
\]
In the case of an isotropic material the strain energy density function $W$ is expressed in terms of the invariants of deformation tensor $I_1$, $I_2$ and $I_3$ \cite{Simon-M,Richter}. 
%defined as
In particular, we will consider the following constitutive law 
\begin{align}\label{used-const.law}
W( F_s)=C_0
+
C_1 \big( I_1-2 \big)
+
C_2 \big( I_2-2 \big)^2,
\end{align}
where $C_0,\ C_1$ and $C_2$ are set as follows
$$
C_0=110 \ \textup{N.cm}^{-2},\ C_1=100 \ \textup{N.cm}^{-2} \ \text{and} \ C_2=110 \ \textup{N.cm}^{-2}.
$$
A pulsatile velocity $ v_\textup {in}$ is enforced on the inlet of the artery. It is given by 
\begin{equation*}
  v_{\text {in}}
 =
\begin{cases}
 5 \ \sin^2(\pi t/0.5) \ \text{cm/s} & \textup{for} \ 5\times(2i) \leq t \leq 5 \times(2i+1),\\
 0 \ \text{cm/s} & \textup{for} \ 5\times(2i+1) \leq t \leq 5 \times(2i+2),
 \end{cases}
 \qquad \textup{for} \ i \in \N^*.
\end{equation*}
It is a periodic continuous function with period $1$ s. It attains its maximum value 5 cm/s at the instants $t=0.25+k$ s, $k \in \N$. 
\subsection{Non-Linear Elastic Modeling of a Stenosed Artery}
\subsubsection*{Blood Flow and Arterial Wall Displacement}
\begin{figure}[h!]
  \begin{subfigure}[b]{0.5\linewidth}
    \centering
    \includegraphics[
   % trim=2cm 1cm 2cm 1cm,height=4.5cm,
    scale=0.3]{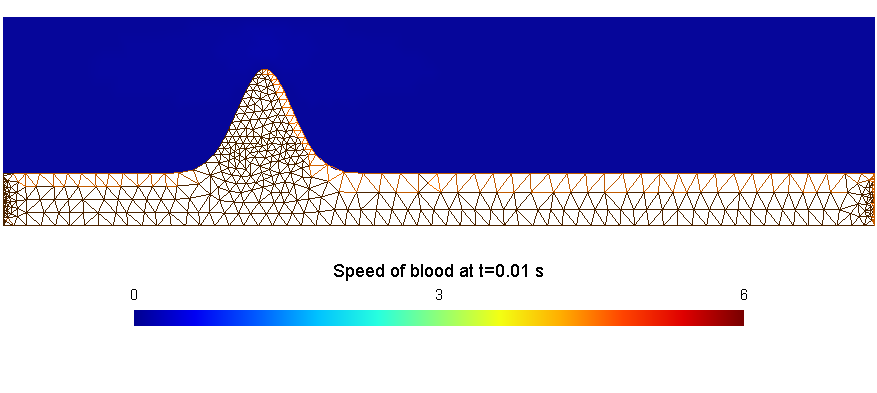}
    \caption{$t=0.01$ s.}
    \label{fig1:a}
    \vspace{1.0ex}
  \end{subfigure}%%
  \begin{subfigure}[b]{0.5\linewidth}
    \centering
     \includegraphics[
   % trim=2cm 1cm 2cm 1cm,height=4.5cm,
    scale=0.3]{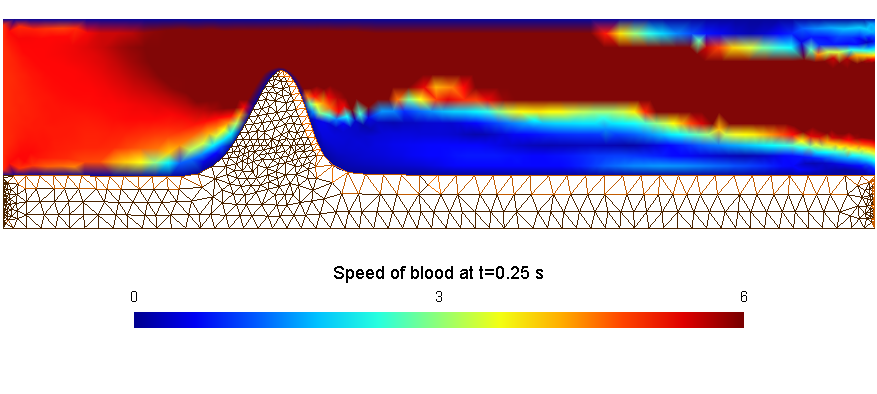}
    \caption{$t=0.25$ s.}
    \label{fig1:b}
    \vspace{1.0ex}
  \end{subfigure}
%\end{figure}
%%% 
%\begin{figure}\ContinuedFloat 
  \begin{subfigure}[b]{0.5\linewidth}
    \centering
     \includegraphics[
   % trim=2cm 1cm 2cm 1cm,height=4.5cm,
    scale=0.3]{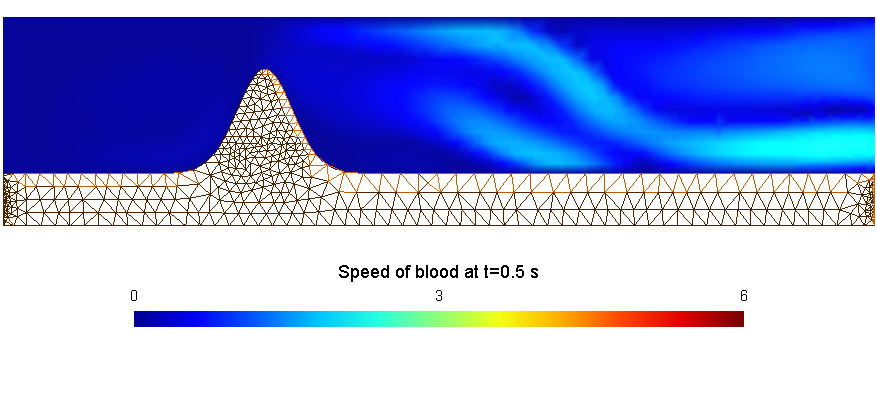}
    \caption{$t=0.5$ s.}
    \label{fig1:c}
  \end{subfigure}%%
  \begin{subfigure}[b]{0.5\linewidth}
    \centering
     \includegraphics[
   % trim=2cm 1cm 2cm 1cm,height=4.5cm,
    scale=0.3]{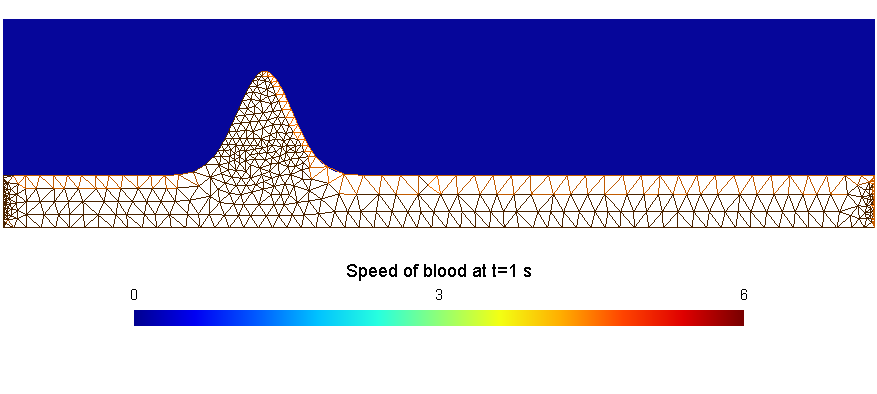}
    \caption{$t=1$ s.}
    \label{fig1:d}
  \end{subfigure}
  \caption{Blood flow in a stenosed artery (cm/s).}
  \label{Flow}
\end{figure}

The first factor that gains our attention is the behavior of the blood flow in the stenosed arteries. Figure \ref{Flow} shows the speed of blood at different instants. One can observe the displacement of the fluid domain, being affected by the displacement of the arterial wall.
\\

\begin{figure}[ht!]
  \begin{subfigure}[b]{0.5\linewidth}
    \centering
    \includegraphics[
%    trim=2cm 1cm 2cm 1cm,height=4.5cm,
    scale=0.3]{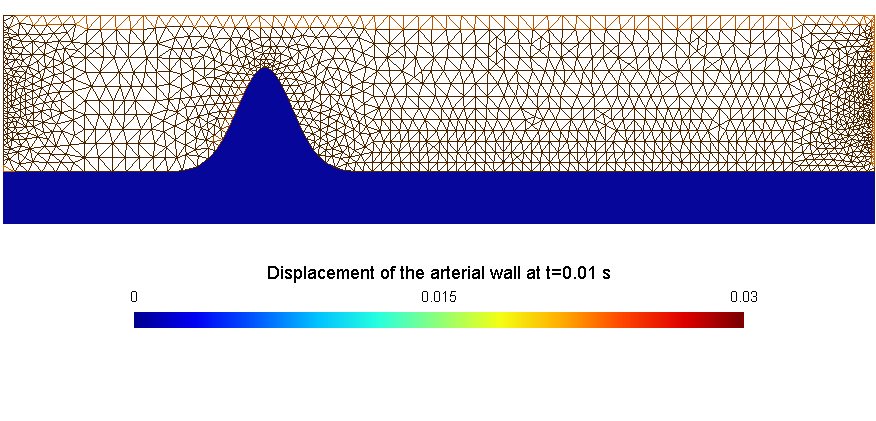}
    \caption{$t=0.01$ s.}
    \label{fig2:a}
    \vspace{1ex}
  \end{subfigure}%%
  \begin{subfigure}[b]{0.5\linewidth}
    \centering
    \includegraphics[
    %trim=2cm 1cm 2cm 1cm,height=4.5cm,
    scale=0.3]{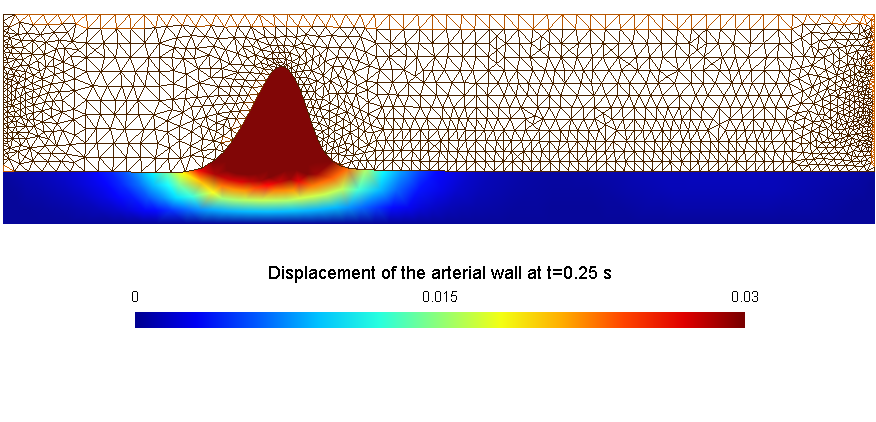}
    \caption{$t=0.25$ s.}
    \label{fig2:b}
    \vspace{1ex}
  \end{subfigure}
%\end{figure}
%%%% 
%\begin{figure}\ContinuedFloat
  \begin{subfigure}[b]{0.5\linewidth}
    \centering
    \includegraphics[
    %trim=2cm 1cm 2cm 1cm,height=4.5cm
    ,scale=0.3]{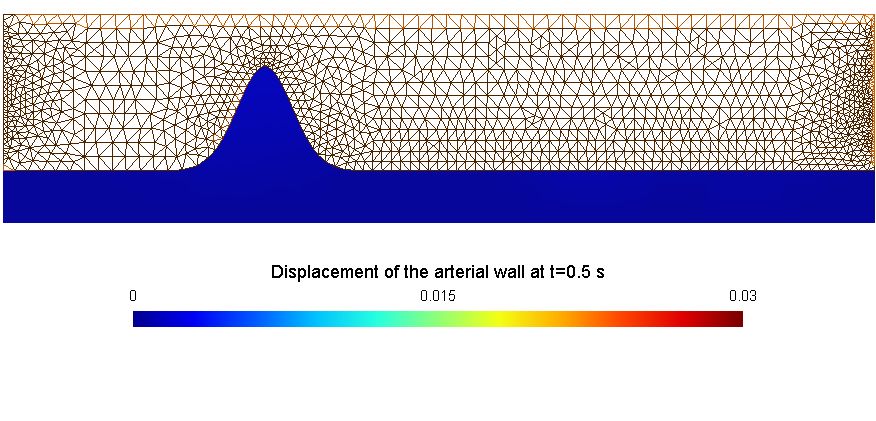}
    \caption{$t=0.5$ s.}
    \label{fig2:c}
  \end{subfigure}%%
  \begin{subfigure}[b]{0.5\linewidth}
    \centering
    \includegraphics[
    %trim=2cm 1cm 2cm 1cm,height=4.5cm,
    scale=0.3]{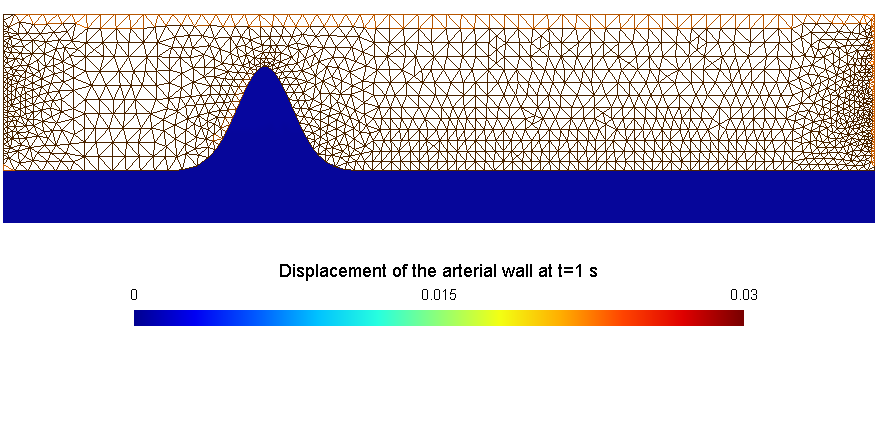}
    \caption{$t=1$ s.}
    \label{fig2:d}
  \end{subfigure}
  \caption{The displacement of the arterial wall (cm).}
  \label{Displacement}
\end{figure}

At the first instant, when the blood starts flowing, its speed is negligible which is observed in Figure \ref{fig1:a}. On the contrast, Figure \ref{fig1:b} shows the remarkable change arising in the lumen of the artery after a time of a quarter of a second. Moreover, the neighborhood of the peak is characterized with a high speed. This is reasonable, indeed, an enforced amount of blood into the artery must pass through it regardless of the diameter of the path or if it is narrowed. Hence, in the narrowed region due to the existence of stenosis, blood speed will become larger. At the instant $t=0.5$, that is when $ v_\textup{in}=0$, the flow through the artery decreases, consequently, the structure domain returns to its equilibrium position as seen in Figure \ref{fig1:c}, however, the effect of the flow is still observed through the lumen domain. This flow behavior continues during the time interval between $0.5$ s and $1$ s to become negligible at $t=1$ s (see Figure \ref{fig1:d}).
Since the displacement of the lumen domain is linked to the displacement of the arterial wall, a view of the displacement of the arterial wall during the same instants will make the observations more obvious. Results during 1 second are given in Figure \ref{Displacement}.
One can observe that the upper part of the stenosis is the region with the highest displacement. In fact, as we reach the peak of the stenosis, i.e, as the stenosis becomes thinner its stiffness will decrease. Consequently, it will be fragile, sensitive to any external force and easily affected by the wall shear stress.
%as well as the speed of blood.
%
%
%

\subsubsection*{Maximum Shear Stress}
\begin{figure}[ht!]
  \begin{subfigure}[b]{0.5\linewidth}
    \centering
    \includegraphics[
    %trim=2cm 1cm 2cm 1cm,height=4.5cm,
    scale=0.3]{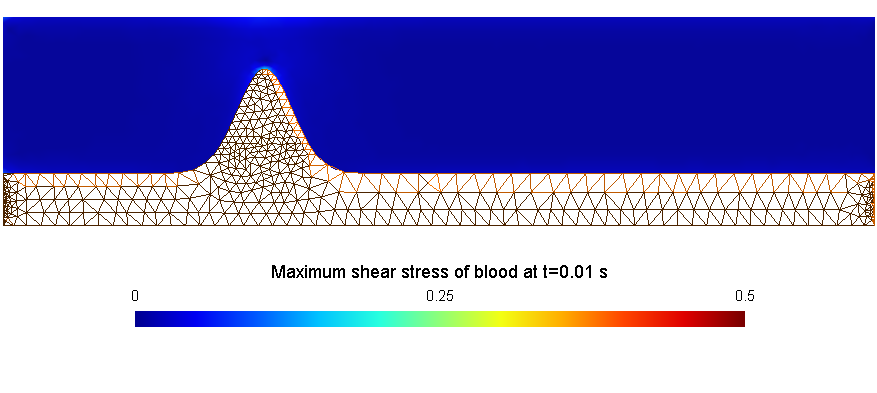}
    \caption{$t=0.01$ s.}
    \label{fig3:a}
    \vspace{1.0ex}
  \end{subfigure}%%
  \begin{subfigure}[b]{0.5\linewidth}
    \centering
    \includegraphics[
   % trim=2cm 1cm 2cm 1cm,height=4.5cm,
    scale=0.3]{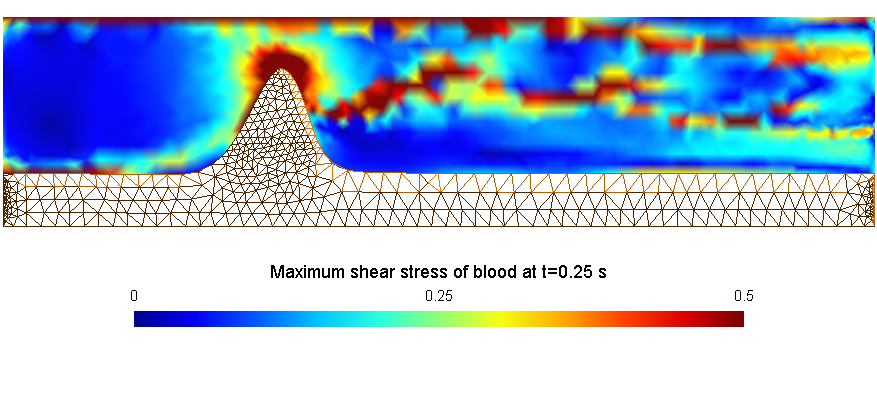}
    \caption{$t=0.25$ s.}
    \label{fig3:b}
    \vspace{1.0ex}
  \end{subfigure}
%\end{figure}
%%%%
%\begin{figure}\ContinuedFloat
  \begin{subfigure}[b]{0.5\linewidth}
    \centering
    \includegraphics[
   % trim=2cm 1cm 2cm 1cm,height=4.5cm,
    scale=0.3]{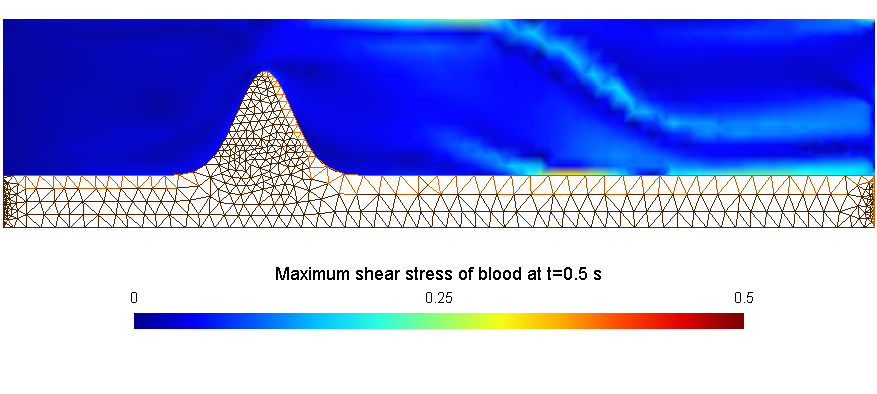}
    \caption{$t=0.5$ s.}
    \label{fig3:c}
  \end{subfigure}%%
  \begin{subfigure}[b]{0.5\linewidth}
    \centering
     \includegraphics[
   % trim=2cm 1cm 2cm 1cm,height=4.5cm,
    scale=0.3]{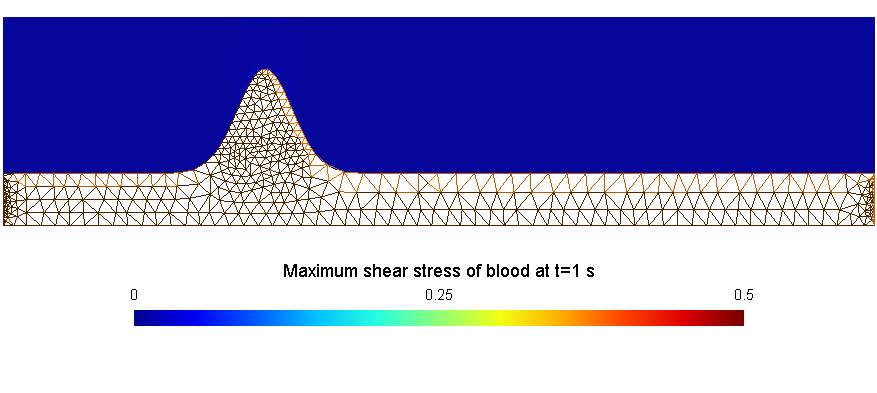}
    \caption{$t=1$ s.}
    \label{fig3:d}
  \end{subfigure}
  \caption{The maximum shear stress (N/cm$^2$).}
  \label{Shear Stress}
\end{figure}

It presents the effect exerted by the fluid on itself. Its expression $\sigma_{max}$ \eqref{max-shs} in terms of the Cauchy stress tensor $ \sigma_f$ reveals its dependence on the strain tensor $ D( v)$ and the pressure of the blood. This means that regions encountering a change in the blood velocity are characterized by a higher maximum shear stress. On the contrary, regions where the values of speed are almost equal are of low maximum shear stress values.  This fact is illustrated on Figure \ref{Shear Stress}.

The maximum shear stress is located in the region of stenosis. More precisely, the observations show that it is located at the peak of the stenosis. In fact, this part of the stenosis is the most fragile part, which makes it more affected by the blood flow.
As the motion of the blood is considered to be periodic; in general it is pulsatilic; and since the structure domain undergoes a deformation, then the stenosis will have an oscillating-like motion. Consequently, a variation in the speed of the blood is recognized. At the instants where the speed in negligible no maximum shear stress is identified (see Figures \ref{fig3:a} and \ref{fig3:d}). On the contrary, at the instant when the stenosis reaches its maximum deformation,  a maximum shear stress of a highest value spotted in the region of stenosis as in Figure \ref{fig3:b}. Further, when the stenosis returns to its equilibrium position, a maximum shear stress is still detected, though, it is of a low value (see Figure \ref{fig3:c}).
\subsubsection{Recirculation Zones}
\begin{figure}[h!]
  \begin{subfigure}[b]{0.5\linewidth}
    \centering
     \includegraphics[
   % trim=2cm 1cm 2cm 1cm,height=4.5cm,
    scale=0.3]{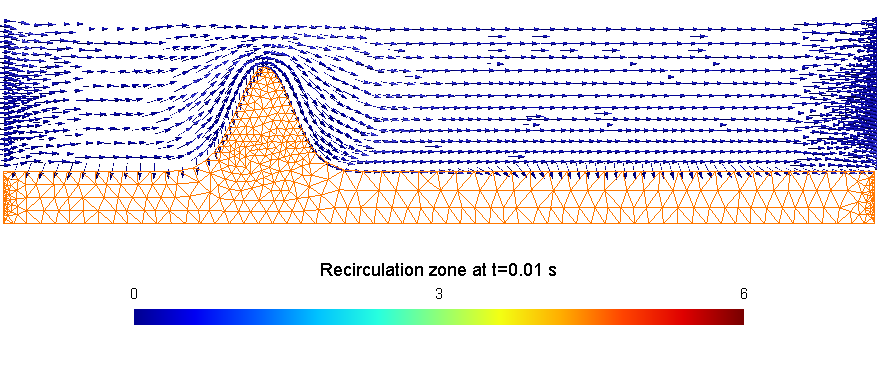}
    \caption{$t=0.01$ s.}
    \label{fig4:a}
    \vspace{1ex}
  \end{subfigure}%%
  \begin{subfigure}[b]{0.5\linewidth}
    \centering
     \includegraphics[
   % trim=2cm 1cm 2cm 1cm,height=4.5cm,
    scale=0.3]{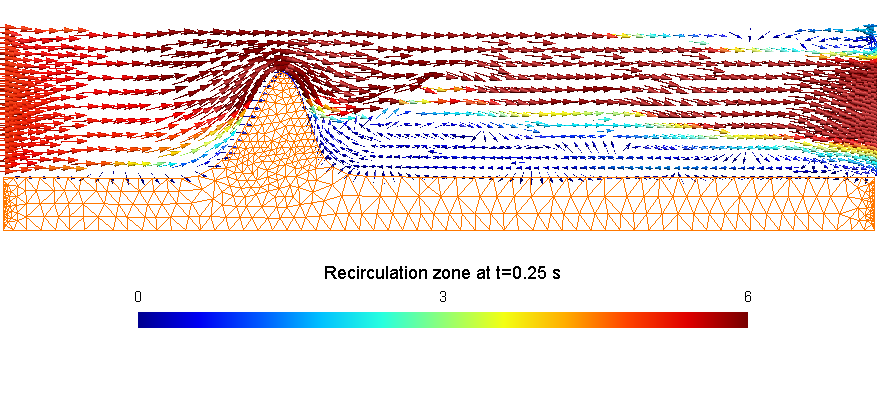}
    \caption{$t=0.25$ s.}
    \label{fig4:b}
    \vspace{1ex}
  \end{subfigure}
%  \end{figure}
%%%%% 
%\begin{figure}\ContinuedFloat
  \begin{subfigure}[b]{0.5\linewidth}
    \centering
     \includegraphics[
   % trim=2cm 1cm 2cm 1cm,height=4.5cm,
    scale=0.3]{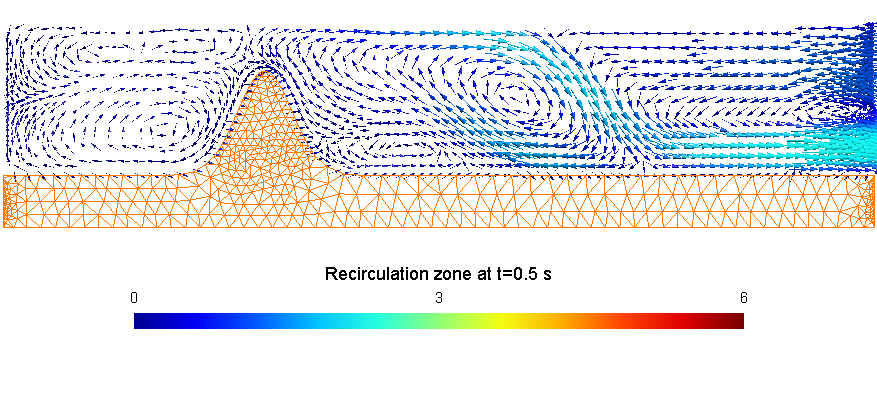}
    \caption{$t=0.5$ s.}
    \label{fig4:c}
  \end{subfigure}%%
  \begin{subfigure}[b]{0.5\linewidth}
    \centering
     \includegraphics[
   % trim=2cm 1cm 2cm 1cm,height=4.5cm,
    scale=0.3]{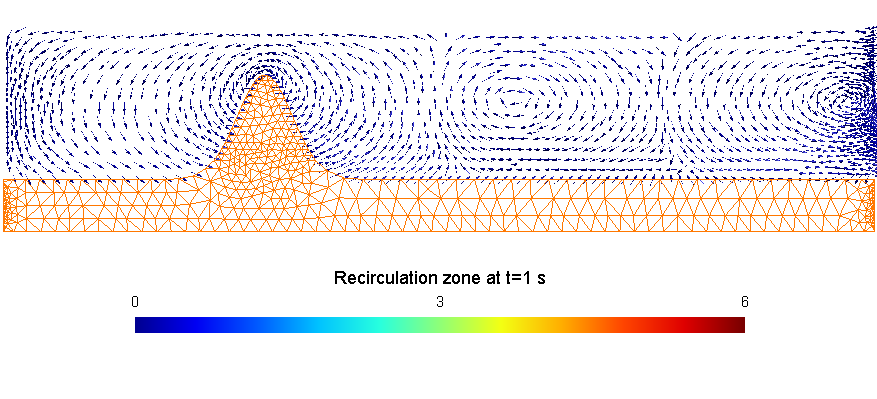}
    \caption{$t=1$ s.}
    \label{fig4:d}
  \end{subfigure}
  \caption{The velocity of blood (cm/s).}
  \label{Rz}
\end{figure}

We are curious about recognizing the recirculation zones. They characterize the regions formed due to the interruption of the flow as a result of the existence of the stenosis and they represent the regions where we have a ripple-like manner.
A vector representation of the blood velocity is illustrated on Figure \ref{Rz} which will help us in configuring the recirculation zones.

Figure \ref{fig4:c} shows the recirculation zones at time $t=0.5$ s. Mainly we can observe a big recirculation zone that is located after the stenosis. The recirculation zone is characterized by a center of negligible speed, which increases as the zone becomes wider. Further, a zone of negligible speed is located between the stenosis and the recirculation zone.
Consequently, in this zone we will observe a phenomenon of sedimentation of the blood. This zone will be the subject of study in Section \ref{solidification-section}. The effect of the recirculation zone on the sedimentation zone would constitute an important tool to build up a rupture model for which the sedimentation zone will be broken due to the forces (arterial wall deformation and blood shear stress) that are applied to it as well as its solid nature.
\\

\begin{figure}[h!]
\centering
\includegraphics[scale=0.4]{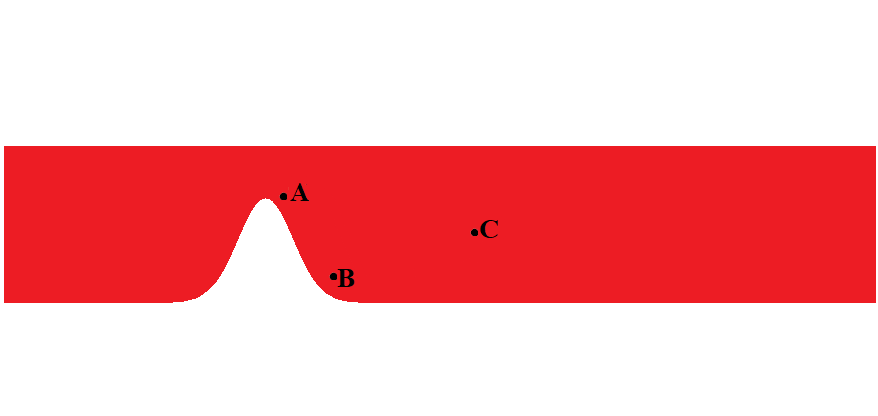}
\caption{Remarkable regions.}
\label{Three zones}
\end{figure}

Previous figures have shown some remarkable regions where the speed varies among them. The same applies for the shear stress and viscosity. In order to study the different phenomena that occur we will focus on three regions. The first region denoted by "A", is located at the peak of the stenosis. Region "B" is located at the adjacent right bottom of the stenosis. And finally, Region "C" is the region including the recirculation zone. 
The three regions are shown on Figure \ref{Three zones}.

\begin{figure}[h!]
    \centering
    \includegraphics[
    %trim=2cm 1cm 2cm 1cm,height=5cm,
    scale=0.4]{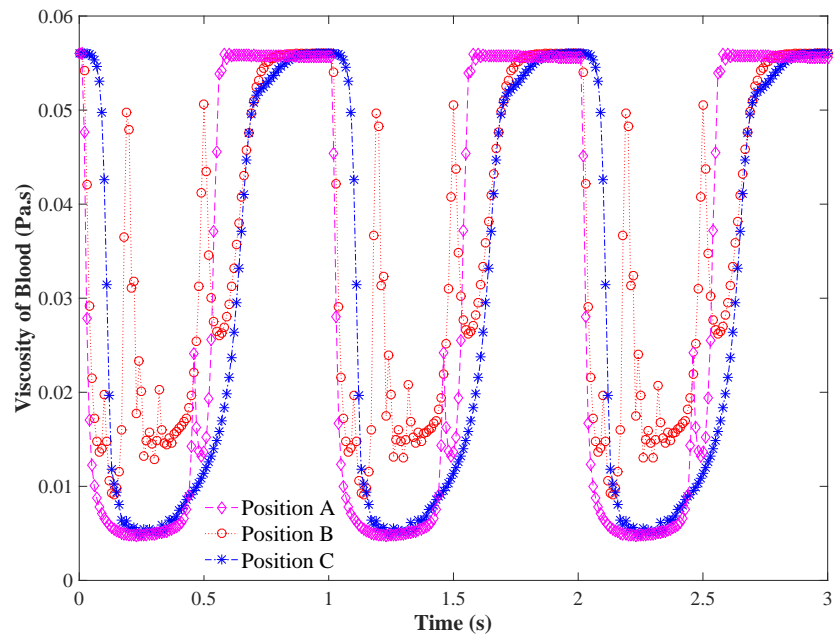}
    \caption{Viscosity of blood.}
    \label{visco-pos1}
%    \vspace{4ex}
  \end{figure}

\begin{figure}[h!]
    \centering
    \includegraphics[
%    trim=2cm 1cm 2cm 1cm,height=5cm,
    scale=0.4]{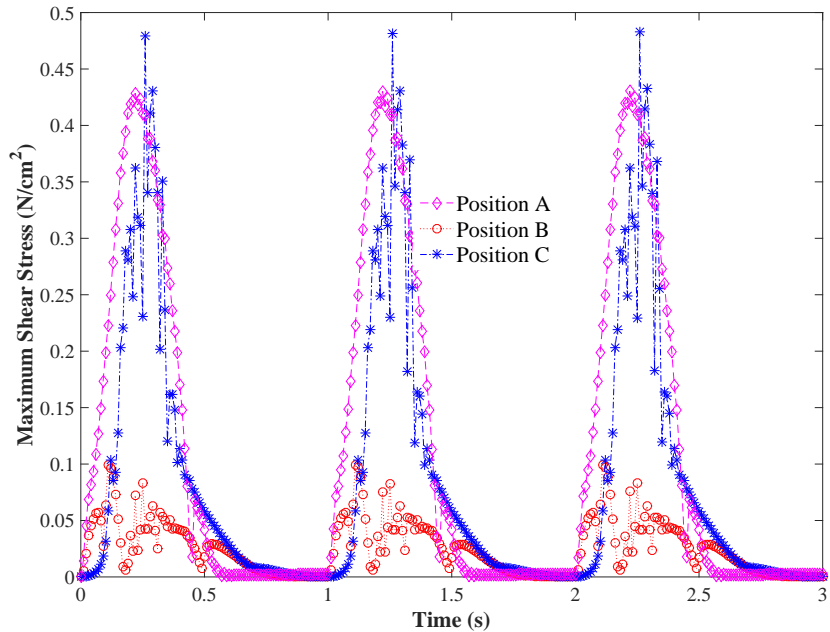}
    \caption{Maximum shear stress.}
    \label{sh-stress1}
    %\vspace{2ex}
  \end{figure}

\begin{figure}[h!]%{0.5\linewidth}
   \centering
   \includegraphics
   [
   %trim=2cm 1cm 4cm 1cm,height=5cm,
   scale=0.4]{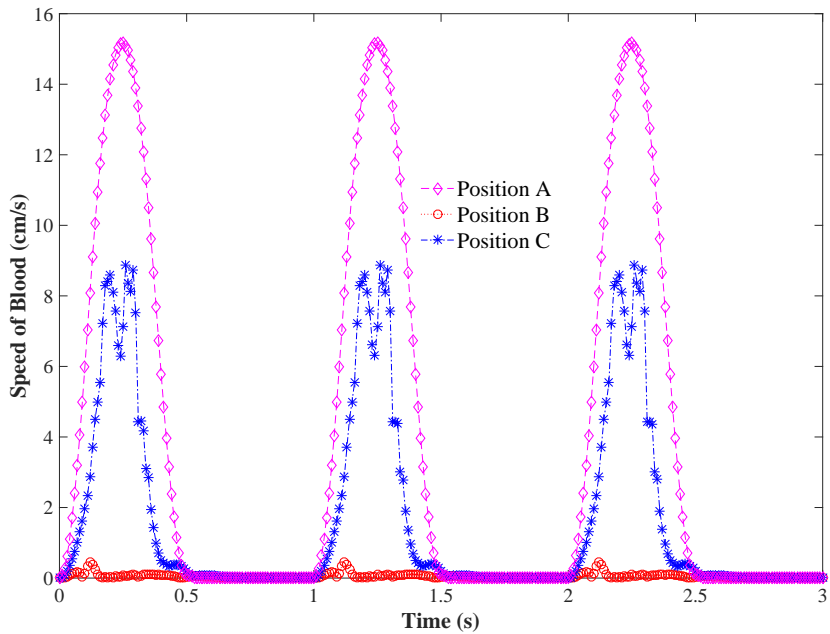}
    \caption{Speed of blood.}
    \label{speed-pos1}
  \end{figure}

The behavior of the viscosity, the shear stress and the speed of the blood at the positions A, B and C are shown on Figures \ref{visco-pos1}, \ref{sh-stress1} and \ref{speed-pos1}, respectively.
%%
%%%%%%%%%%%%%%%%%%%%%%%%%%%%%%%%%%%%%%%%%%%%%%%%%
 %%
\newpage
Figures \ref{visco-pos1}, \ref{sh-stress1} and \ref{speed-pos1} show that the position B is characterized by a high viscosity and a negligible speed, which will lead to the formation of a more viscous region. Consequently, the shear stress at this position is small. In the next section, this region will be identified as a solidification zone. Whereas, the position A is characterized by a high speed of blood with a moderate maximum shear stress rate. In fact, due to the narrowing of the artery at the region of stenosis, the blood speed will be high, as we mentioned previously. On the other hand, at the position C, located in the recirculation zone, we have a high shear stress, due to the change encountering on the speed in this region.
The investigation of these variables at these remarkable positions will help us in locating the solidification zone, where a clot would form, and analyze the external forces and stresses applied to it.

\section{Solidification of Blood and its Rupture}
\label{solidification-section}
Blood clots are formed whenever the flowing blood come in contact with a foreign substance in the skin or in the blood vessels wall. 
They can be classified into two types: thrombosis, which are stationary clots, though they can cause the blockage of a flow; embolisms, which detach into the blood flow and can, somewhere in a site far away from the thrombosis, block the flow. 
This type of clots is dangerous and causes infarctions, more precisely, if the blockage occurs in the brain it results a stroke, if it occurs in the heart a heart attack would result, or in the lungs it would cause a pulmonary embolism. 
In particular, in the situations where plaques formed from fats, lipids, cholesterol or other foreign substances found in the blood are identified, over time, they harden causing the narrowing of the artery.

\subsection{Detection of the solidification Zone} 
The final step in the formation of plaque- which is rupture- does not always occur. We believe that it is linked to the solidified blood and is influenced by some factors that we will discuss later from the numerical viewpoint. For this reason, the first step in building a rupture model is characterized by spotting the region of solidification. That is, we will investigate based on the rheology of blood, constitutive models and the numerical results of Section \ref{Numerical Results}, the region where the blood transits into a gel state. In general, in vivo, blood is liquid in state, then a change in its state from liquid to gel is linked to a change in the viscosity. Indeed, as the viscosity increases, a more solidified material is acquired. Hence, a solidification zone should be identified by a sufficiently large viscosity. Though, we can detect many regions that are characterized by high viscosity values. In fact, in regions where the values of the velocities are almost equal, the viscosity is of high values, this fact arises from the relation between the viscosity $\mu$ and the deformation tensor $ D( v)$. In other words, as the rate of change of the velocity expressed by $ D( v)$ is negligible then the viscosity tends to reach its highest asymptotic value $\mu_0$.
This reveals that the condition of possessing a high viscosity is insufficient to detect the solidification zone. Hence, another condition is essential to achieve a precise location of the solidification zone. It is recognized that gel and jelly-like materials spread and flow slowly. Consequently, the solidification zone must also obey the fact that it is of negligible speed. 

Results of Section \ref{Numerical Results} have shown the existence of recirculation zones after the stenosis due to the blockage of flow by the stenosis. These zones are characterized by a negligible flow resulting from a negligible blood velocity at its center, which increases as the circles formed in this zone become larger in diameter. Between the recirculation zone and the stenosis, in particular, at the edge of stenosis, we detect a region where the flow is of negligible speed and of a high viscosity. This region is identified as the solidification zone since it possesses the characteristics mentioned above. 

Notice that, as the formation of atherosclerosis is a long-time process, the formation of the solidification region is as well. Literally, viscosity is a time-dependent function. In fact, the formation of these regions depends on the flow of blood at each pulse, that is, the viscosity depends on its history. 

Numerically, we consider a threshold value $\mu_{th}$ of the viscosity such that when the computed viscosity $\mu$ exceeds $\mu_{th}$ a region $\mathcal{D}_\mu$ of high viscosity is identified. Similarly, we consider a threshold value $v_{th}$ of the blood speed. If the speed $\norm  v \norm_{2}$ is such that it is less than $v_{th}$ then we locate a region $\md_v$ of a negligible velocity. The solidification region $\mathcal{R}_{\mathfrak{s}}$ is the intersection of the two located regions $\md_\mu$ and $\md_v$. More precisely, the region $\mathcal{R}_{\mathfrak{s}}$ satisfies possessing a high viscosity as well as a negligible speed.
\\

We model the blood using the Carreau model. We recall that its associated viscosity is given by \eqref{Carreau}.

A modification is applied to Carreau model so that the viscosity becomes a time-dependent function expressed in terms of its history. As we are performing iterative simulations, then at each time iteration $k$, for $k \in \N$, we will express the viscosity $\mu_k$ in terms of the local in time average viscosity $\hat{\mu}_{k}$ given by
\begin{equation*}
\hat{\mu}_k =
\begin{cases}
\mu_0 =0.056  &\qquad \textup{if} \ k=0,\\
0.035 		  &\qquad \textup{for} \ 1\leq k \leq 4,\\
\displaystyle \dfrac{1}{5} \sum_{i=1}^5 \mu_{k-i} &\qquad \textup{for} \ k \geq 5,
\end{cases} 
\end{equation*}
where $\mu_{k-i}$ represents the viscosity of blood at each iteration $k-i$. 
Then at each time iteration $k \in \N$, we set
\[
\mu_{\infty,k}= \mu_\infty - t \times 10^{-3} \times \hat{\mu}_{k}^{0.2}
\]
and 
\[
\mu_{0,k}= \mu_0 - t \times 10^{-3} \times \hat{\mu}_{k}^{0.2},
\]
where $t= k \times \Delta t$ with $\Delta t = 10^{-2}$ s is the time step.
\\
Hence, at the iteration $k$, the viscosity expression becomes
\begin{align}\label{Modified Carreau model}
\mu_{k}=
\mu_{\infty,k} +(\mu_{0,k} - \mu_{\infty,k})
\big[1+ (\lambda \ocg_k)^2 \big]^{\frac{n-1}{2}},
\end{align}
with $\ocg_k$ given by \eqref{shear rate} as 
$
\sqrt{2
\textup{tr}
\big( D( v^{k-1})
\big)^2}$.
In a two dimensional space, its explicit expression is
\[
\sqrt{
	2
	\bigg(
	\dfrac{d}{dx} v^{k-1}_{1}
	\bigg)^2
   +2
   \bigg(
   \dfrac{d}{dy} v^{k-1}_2
   \bigg
   )^2
   +
   \bigg(
    \dfrac{d}{dy} v^{k-1}_1 + \dfrac{d}{dx} v^{k-1}_2
    \bigg)^2
  	}.
\]
It should be noticed from the context that the superindices $k-1$ and $k$ refer to the time iteration, while the subindices $1$ and $2$ stand for the vector components of $ v$.
\\

\begin{figure}[h!]
\centering
\begin{subfigure}{.5\textwidth}
  \centering
  \includegraphics[scale=0.7,width=1.\linewidth]{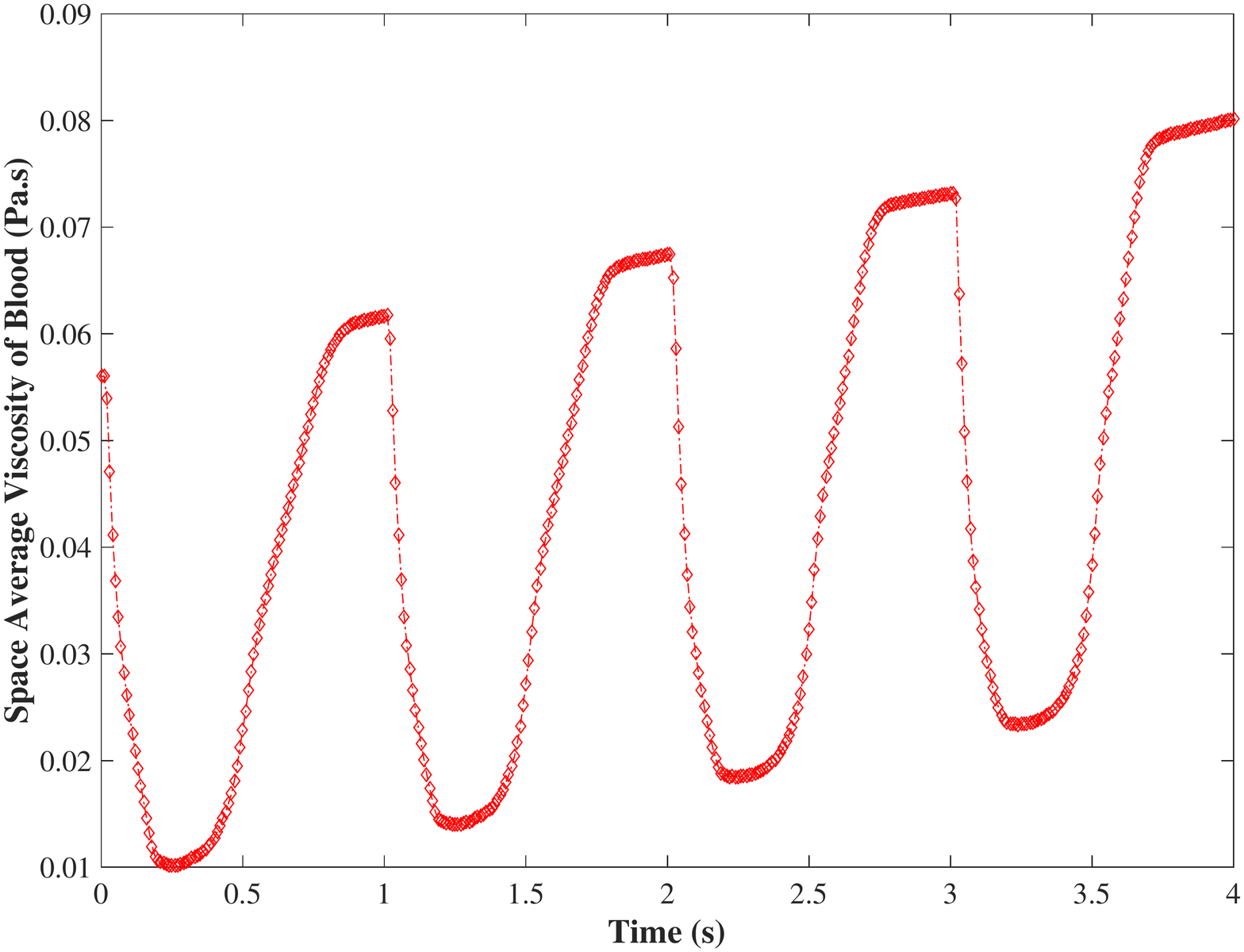}
  \caption{Space average viscosity of blood.}
  \label{viscosity_nwtn and non_nwtn}
\end{subfigure}%
\begin{subfigure}{.5\textwidth}
  \centering
  \includegraphics[scale=0.7,width=1.\linewidth]{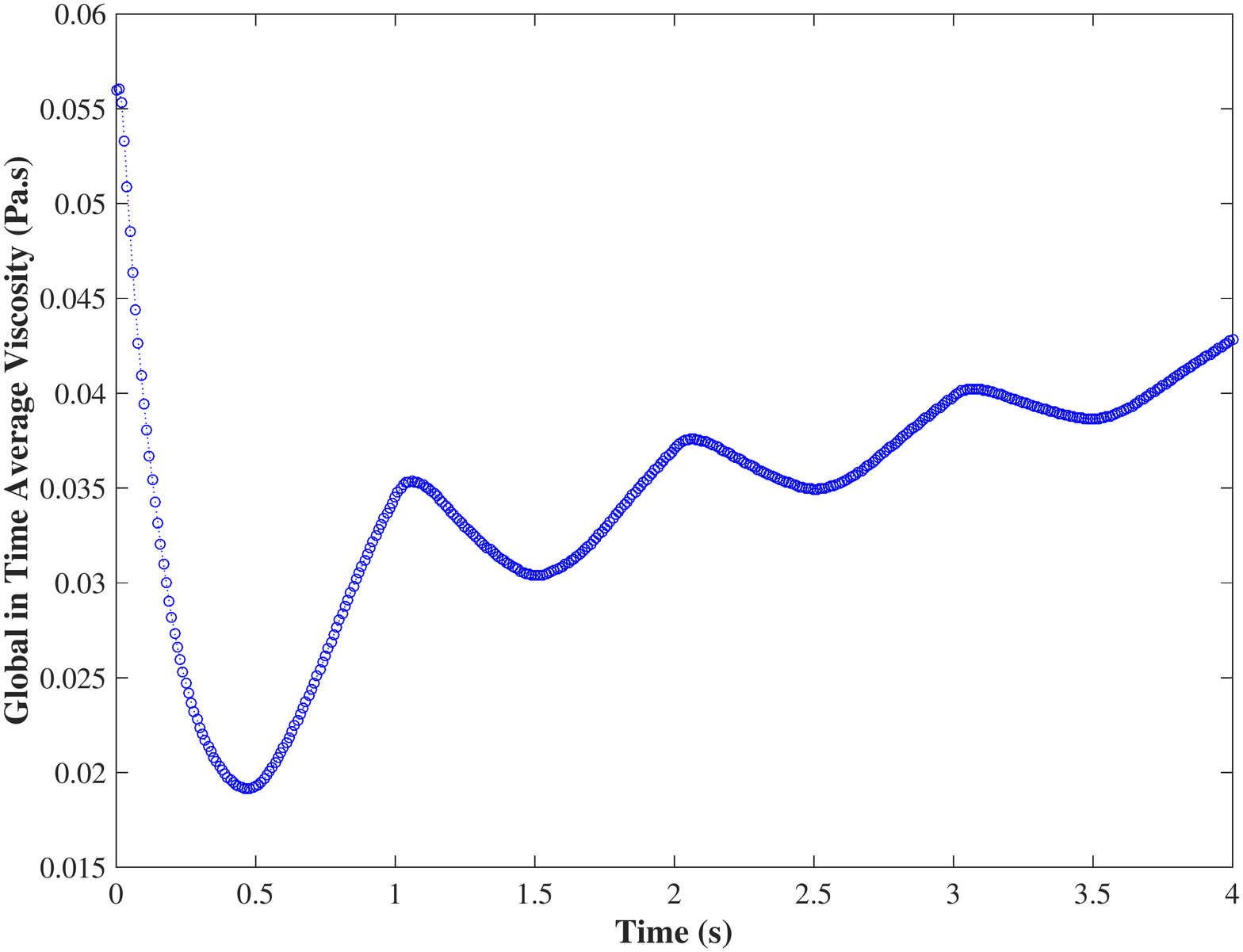}
  \caption{Global in time average viscosity of blood.}
  \label{av. viscosity_nwtn and non_nwtn}
\end{subfigure}
\caption{Space average viscosity (left) and global in time average viscosity (right) of a non-Newtonian blood.}
\label{visc and av visc}
\end{figure}

To set the threshold values $\mu_{th}$ and $v_{th}$, we plot the most remarkable values on a specified time interval $[t_0,T]$, $t_0>0$. The highest remarkable value is set to be the threshold in case of viscosity. Whereas, in case of speed we consider the lowest remarkable value.
At each iteration $k$, the global in time average viscosity $\overline{\mu}_k$ is given by the following relation
\begin{align*}
\overline{\mu}_k= \dfrac{1}{k+1} \sum_{i=0}^k \mu_{i},
\end{align*}
where $\mu_i$ represents the viscosity of the blood at iteration $i$, $0 \leq i \leq k$ and we set $\overline{\mu}_0=\mu_0 = 0.056$ Pa.s.
The graphs corresponding to the viscosity and the average viscosity obeying \eqref{Modified Carreau model} during 3 seconds are plotted on Figure \ref{visc and av visc}.
\\

\begin{figure}[h!]%{0.5\linewidth}
   \centering
   \includegraphics
   [
   %trim=3cm 0cm 6cm 1cm,height=6cm,
   scale=0.4]{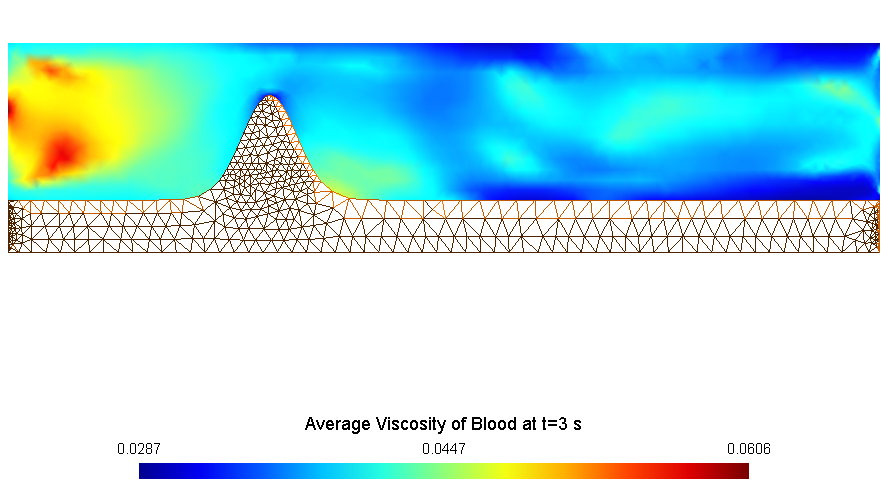}
    \caption{Average viscosity of blood at time $t_0=3$ s.}
    \label{Av_visco_carreau}
  \end{figure}

\begin{figure}[h!]%{0.5\linewidth}
   \centering
   \includegraphics
   [
   %trim=2cm 1cm 4cm 1cm,height=5cm,
   scale=0.4]{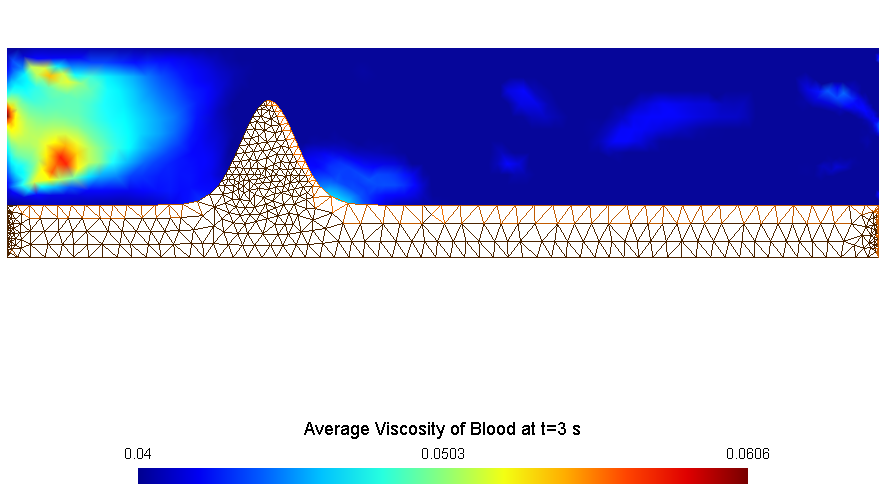}
    \caption{Regions of average viscosity greater than 0.04 Pa.s.}
    \label{rescaling viscosity}
  \end{figure}

\begin{figure}[h!]%{0.5\linewidth}
   \centering
   \includegraphics
   [
   %trim=2cm 1cm 4cm 1cm,height=5cm,
   scale=0.4]{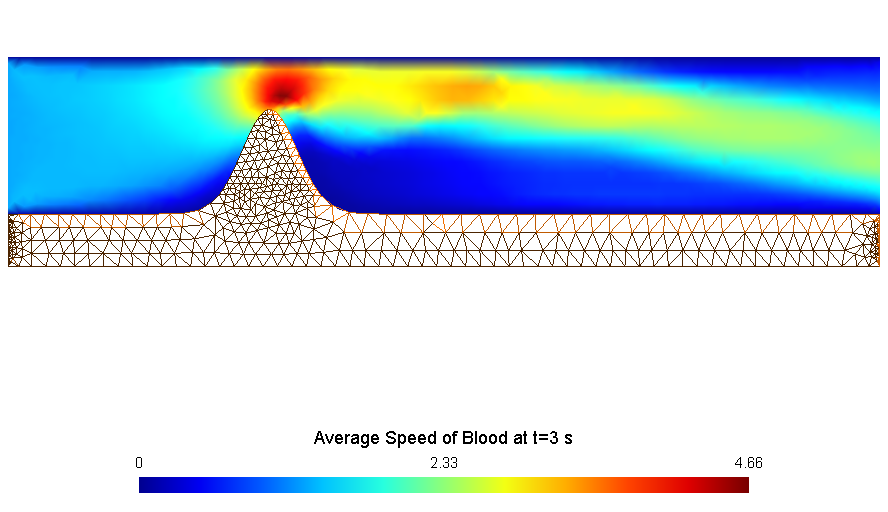}
    \caption{Average speed of blood at time $t_0=3$ s.}
    \label{Av-speed-it300}
  \end{figure}

The pattern of the average viscosity at the instant $t_0=3$ s is illustrated on Figure \ref{Av_visco_carreau}.
 %%
From Figure \ref{Av_visco_carreau} we observe mainly two regions possessing high average viscosity. The first region, located at the inlet of the artery, is a region where the particles constituting it are characterized by values of the velocity that are almost equal. Thus, based on the expression of the viscosity expressed in terms of the deformation tensor $ D( v)$ a high viscosity results. On the other hand, the second region of a high average viscosity, is located near the stenosis. The existence of stenosis prevents the flow from reaching the spot at the edge of the stenosis. Consequently, blood will become more viscous. In particular, the values of the viscosity in these two remarkable regions are greater than $0.04$ Pa.s. Rescaling the data we get a precise location of the regions which are characterized by an average viscosity greater than $0.04$ Pa.s (see Figure \ref{rescaling viscosity}).
 %%
As a result, we set the threshold $\mu_{th}$ to be $0.04$ Pa.s.
As we mentioned previously, another condition is desired to obtain a precise location of the solidification zone. More precisely, the zone must be characterized by a negligible speed. The average speed at time $t_0=3$ s is illustrated on Figure \ref{Av-speed-it300}.
 %%
\\

\begin{figure}[h!]%{0.5\linewidth}
   \centering
   \includegraphics
   [
   %trim=2cm 1cm 4cm 1cm,height=5cm,
   scale=0.4]{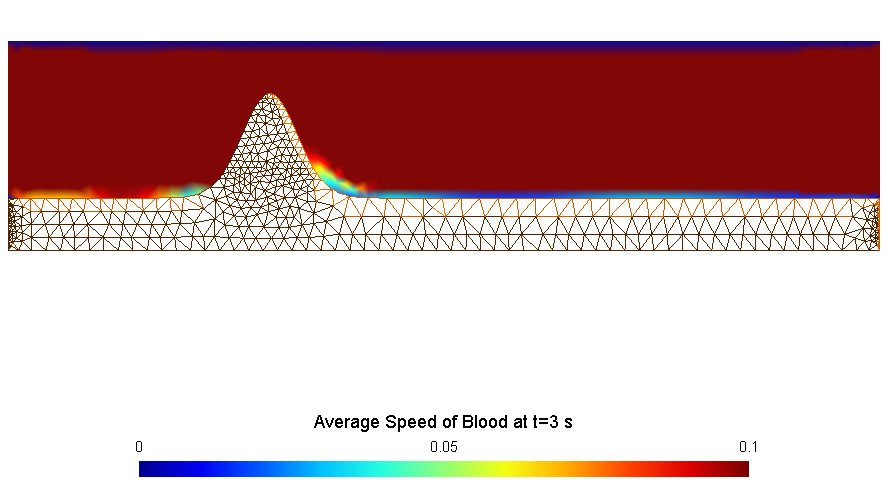}
    \caption{Regions of average speed less than 0.1 cm/s.}
    \label{rescaling-Av-speed-it300}
  \end{figure}

\begin{figure}[h!]
        \centering
        \begin{tikzpicture}
            [path image/.style={path picture={\node at (path picture bounding box.center) {\includegraphics[height=3cm]{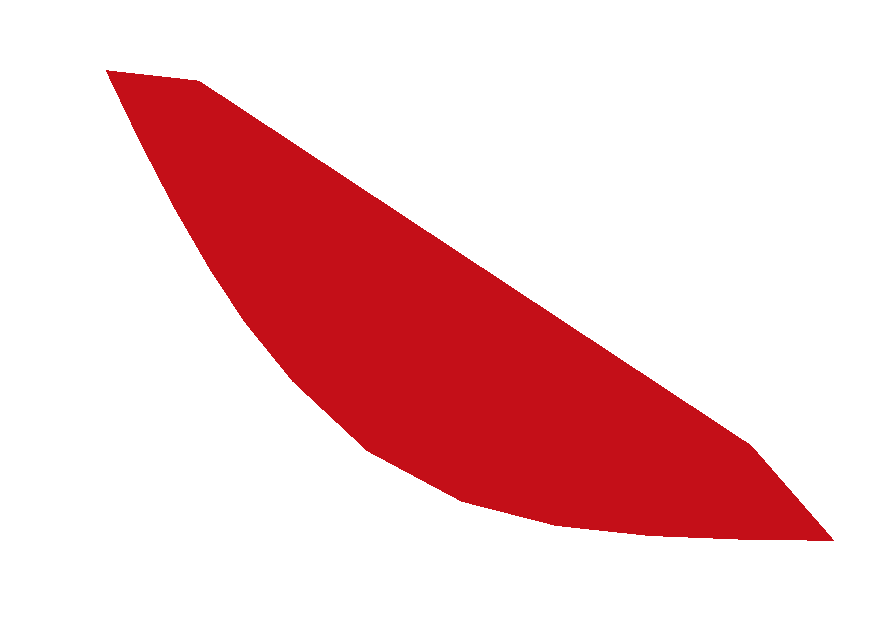}};}}]

            \node (img) {\includegraphics[width=.7\linewidth]{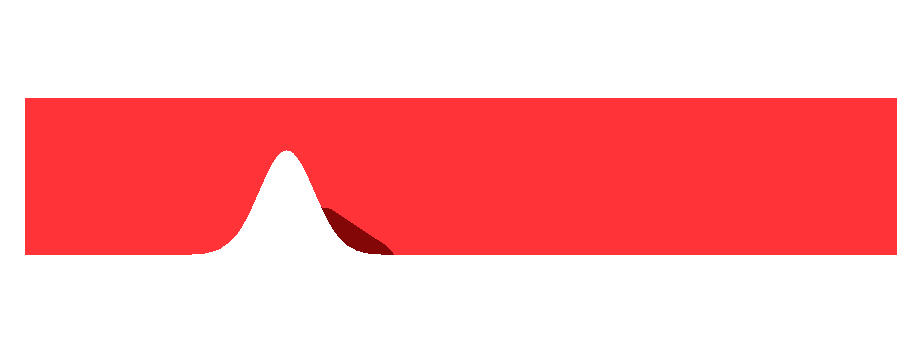}};
            \node (c1) [draw, circle, black, text width=1.cm] at (img.center) at (-1.22,-0.7) {};
            \draw [black] (c1.east) -- (img.east);
            \draw [path image=example-image-a,draw=black,thick] (img.east) circle (2cm);
        \end{tikzpicture}
        \caption[The solidification zone $\mathcal{R}_{\mathfrak{s}}(t)$.]{The solidification zone $\mathcal{R}_{\mathfrak{s}}(t)$.}
        \label{CoagZone}
    \end{figure}
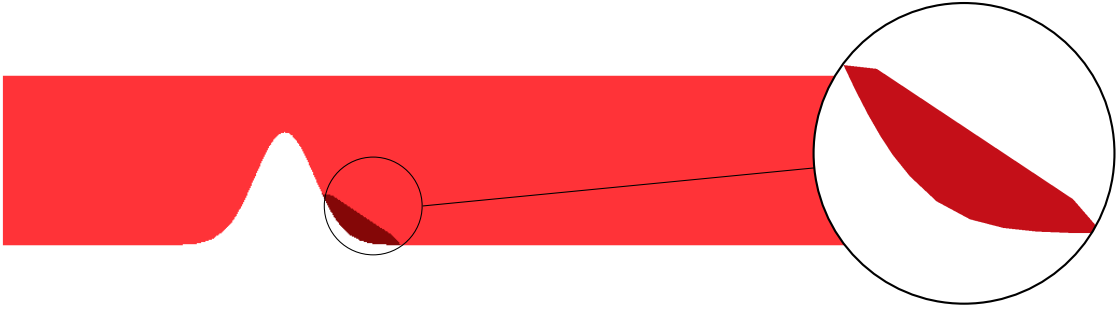

Figure \ref{Av-speed-it300} shows that the average speed attains its highest value above the peak of the stenosis. In contrast, the lowest value is configured at the edge of the stenosis. A rescaling of the values would help us get a precise data. Indeed, Figure \ref{rescaling-Av-speed-it300} reveals that the speed of the blood existing at the edge of the stenosis is of maximum value 0.1 cm/s.
 %%
To sum up, at $t_0=3$ s, Figures \ref{rescaling viscosity} and \ref{rescaling-Av-speed-it300} showed a region at the edge of the stenosis where the blood is characterized by a high viscosity and a low speed. In fact, by setting $\mu_{th}=0.04$ Pa.s and $v_{th}= 0.1$ cm/s we get an accurate detection of the solidification zone at the edge of the stenosis. 
The solidification zone is given in Figure \ref{CoagZone}.

%
%\begin{figure}[h!]
%  \centering
%  \includegraphics[scale=0.5,width=.5\linewidth]{Cz.png}
%\caption{The solidification zone $\mathcal{R}_{\mathfrak{s}}(t)$.}
%\label{CoagZone}
%\end{figure}

\subsection{Forces Acting on the solidification Zone}

\begin{figure}[h!]%{0.5\linewidth}
   \centering
   \includegraphics
   [
   %trim=2cm 1cm 4cm 1cm,height=5cm,
   scale=0.35]{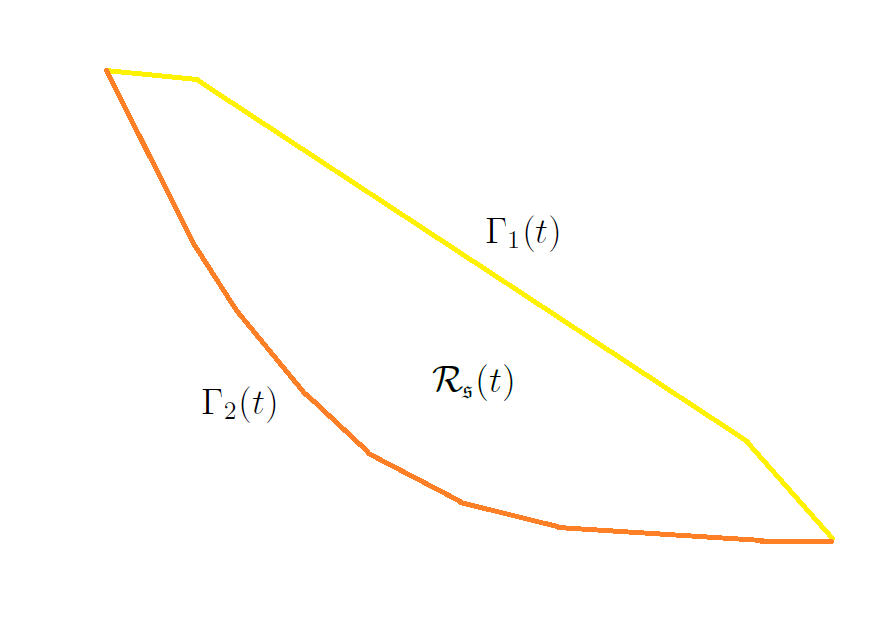}
    \caption{The domain of the solidification zone.}
    \label{domain-coag}
  \end{figure}

Having located the region of solidification $\mathcal{R}_{\mathfrak{s}}(t)$, $t>t_0>0$ , we proceed to identify the factors that will lead to the rupture of the semi-solidified blood. The solidification region is made up of blood in gel state, which, similarly as the plaque will be under the effect of a force exerted by the pressure of the blood and the shear stress. Moreover, being located at the edge of the stenosis, it will be affected by the stenosis displacement.
At any instant $t>0$, we denote by $\Omega_f(t)$ the domain corresponding to the lumen of the artery and by $\Omega_s(t)$ the domain representing the arterial wall. On $\Omega_f(t)$ we define $ v$ the velocity of the blood and $p_f$ to be its pressure. On the other hand, the motion of the arterial wall is defined by its displacement $ \xi_s$. The boundary $\partial \mathcal{R}_{\mathfrak{s}}(t)$ of the  solidification zone $\mathcal{R}_{\mathfrak{s}}(t)$ is decomposed into $\Gamma_1(t)$ and $\Gamma_2(t)$ as shown on Figure \ref{domain-coag}.  
 %%
\subsubsection*{Linear Elasticity of the Semi-Solidified Blood}
We deal with the solidification zone from the perspective of being an elastic material that obeys Hooke's law. Let us designate by $ u=(u_1,u_2)$ the displacement of the domain $\mathcal{R}_{\mathfrak{s}}(t)$. The solidification zone is under the effect of an external surface force. In particular, a surface force $f_c$, representing the shear stress is applied from the blood surrounding the solidification zone to the boundary $\Gamma_1(t)$. Thus, the expression of $f_c$ is given in terms of the Cauchy stress tensor $ \sigma_f( v, p_f)$ by
\begin{align}\label{Surface force_on solidification zone}
 f_c = -  \sigma_f( v,p_f) \  n_f  \qquad \textup{on} \quad \Gamma_1(t) \times (t_0,T). 
\end{align}
where $ n_f$ is the outward normal to the domain $\Omega_f(t) \setminus \mathcal{R}_{\mathfrak{s}}(t)$ and $\sigma_f(v,p_f)$ is given by \eqref{Cauchy-stress}.
\\
On the other hand, since the border $\Gamma_2(t)$ constitutes a part of the common boundary $\Gamma_c(t)=\partial \Omega_f(t) \cap \partial \Omega_s(t)$ then we must ensure the continuity of the deformation on this boundary, that is, we impose the condition
\begin{align}\label{Continuity of Displacement}
 u =  \xi_s \qquad \textup{on} \quad \Gamma_2(t) \times (t_0,T).
\end{align}
 
As a result, the elasticity equations satisfied by the displacement  $u$ of the solidification zone $\mathcal{R}_{\mathfrak{s}}(t)$ are
%%%
\begin{equation}\label{Elasticity eqn_CoagZone}
\begin{cases}

- \textbf{div} \   \sigma_c( u) = 0 & \qquad \textup{in} \quad \mathcal{R}_{\mathfrak{s}}(t) \times (t_0,T), \\

 \sigma_c( u) \  n_c =  f_c 
& \qquad \textup{on} \quad \Gamma_1(t) \times (t_0,T), \\

 u =  \xi_s & \qquad \textup{on} \quad \Gamma_2(t) \times (t_0,T),
\end{cases}
\end{equation}
%%%
where $ n_c$ is the outward normal to the solidification zone $\mathcal{R}_{\mathfrak{s}}(t)$. The Cauchy stress tensor $ \sigma_c( u)$ is expressed in terms of the strain tensor $ \varepsilon( u)=\dfrac{1}{2}\big( \nabla  u + ( \nabla  u)^t \big)$ by Hooke's law
%%%
\begin{align}\label{Hooke's law}
 \sigma_c( u)
=
2 \mu_c  \varepsilon( u)+\lambda_c \textup{tr}( \varepsilon( u)) \ \Id
\end{align}
with $\mu_c$ and $\lambda_c$ are the Lam\'e constants that are given in terms of the Young's modulus $E$ and the Poisson's ratio $\nu$ as
\begin{align*}
\lambda_c = \dfrac{\nu E}{(1-2\nu)(1+\nu)} \qquad \textup{and} \qquad \mu_c=\dfrac{E}{2(1+\nu)}.
\end{align*}

For a clot, which is is assumed to be an incompressible material, the Poisson's ratio is $\nu=0.492$ \cite{Adam}. Further, its Young's modulus (Elastic modulus) $E =14.5$ MPa \cite{Collet2005}.  
\\

\begin{figure}[h!]
\centering
\begin{subfigure}{.5\textwidth}
  \centering
  \includegraphics[scale=0.7,width=1.\linewidth]{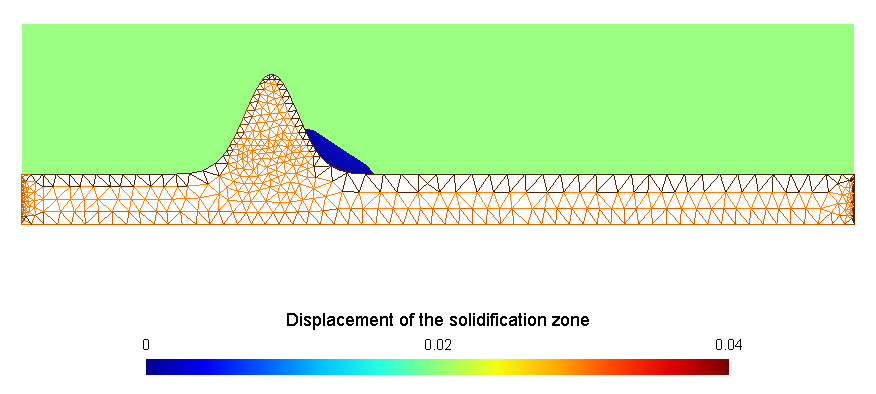}
  \caption{The displacement of the solidification zone at $t=3$ s.}
  %\label{viscosity_nwtn and non_nwtn}
\end{subfigure}%
\hfill
\begin{subfigure}{.5\textwidth}
  \centering
  \includegraphics[scale=0.7,width=1.\linewidth]{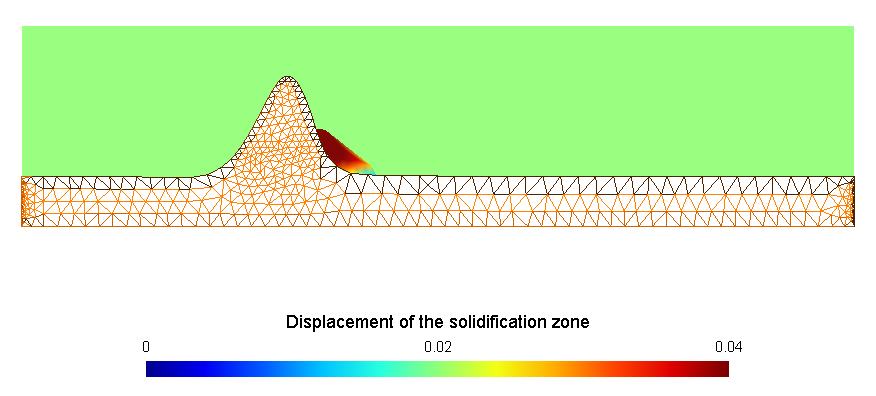}
  \caption{The displacement of the solidification zone at $t=3.25$ s.}
  %\label{av. viscosity_nwtn and non_nwtn}
\end{subfigure}
\caption{The deformation of the solidification zone between $t=3$ s and $t=3.25$ s.}
\label{disp_coag}%
\end{figure}

In order to write the variational formulation associated to System \eqref{Elasticity eqn_CoagZone}, we rewrite it as a partial differential equation with homogeneous Dirichlet boundary condition. For this reason, we consider a function $\mathfrak{h} \in H^1(\mathcal{R}_{\mathfrak{s}}(t))$ such that $\gamma_0 (\mathfrak{h})=  \xi_s$, where 
\begin{align*}
\gamma_0 : H^{1/2}(\Gamma_2(t)) \mapsto H^1(\mathcal{R}_{\mathfrak{s}}(t))
\end{align*}
is the trace operator.
\\
Take $ \zeta= u-\mathfrak{h}$ which is a function in $H^1(\mathcal{R}_{\mathfrak{s}}(t))$ that vanishes on $\Gamma_2(t)$. Since $ \sigma_c( u)$ is a function of $ \varepsilon( u)$ which is linear, then we have
\begin{align*}
 \sigma_c( \zeta)=
 \sigma_c( u)-
 \sigma_c(\mathfrak{h}).
\end{align*}
Therefore, System \eqref{Elasticity eqn_CoagZone} is equivalent to
\begin{equation}\label{Elasticity eqn_CoagZone-Homo}
\begin{cases}

- \textbf{div} \   \sigma_c( \zeta) = \textbf{div} \   \sigma_c(\mathfrak{h}) & \qquad \textup{in} \quad \mathcal{R}_{\mathfrak{s}}(t) \times (t_0,T), \\

 \sigma_c( \zeta) \  n_c =  f_c -   \sigma_c(\mathfrak{h}) \  n_c
& \qquad \textup{on} \quad \Gamma_1(t) \times (t_0,T), \\

 \zeta =  0 & \qquad \textup{on} \quad \Gamma_2(t) \times (t_0,T).
\end{cases}
\end{equation}
\\
The variational formulation associated to System \eqref{Elasticity eqn_CoagZone-Homo} is derived by considering a test function 
\[
 \eta_c \in \mathcal{W}_c =
\{  \eta \in H^1(\mathcal{R}_{\mathfrak{s}}(t)), \  \eta =  0 \ \textup{on} \ \Gamma_2(t) 
\}
\]
to get
\begin{align}\label{V.F_Coag-zone}
\int_{\mathcal{R}_{\mathfrak{s}}(t)}
 \sigma_c( \zeta) :  \nabla  \eta_c \ d x
-
\int_{\Gamma_1(t)} 
 \sigma_c( \zeta) \  n_c \cdot  \eta_c \ d\Gamma
= 
\int_{\mathcal{R}_{\mathfrak{s}}(t)}
[\textbf{div} \  \sigma_c(\mathfrak{h})] \cdot  \eta_c \ d x
.
\end{align}
Substituting $ \sigma_c( \zeta)$ by its expression \eqref{Hooke's law} and $ f_c$ by \eqref{Surface force_on solidification zone} we can rewrite \eqref{V.F_Coag-zone} as
\begin{equation}\label{V.F-elasticity eqn}
\begin{cases}
\displaystyle
2 \mu_c 
\int_{\mathcal{R}_{\mathfrak{s}}(t)}
 \varepsilon( \zeta):  \varepsilon( \eta_c) \ d x
+
\lambda_c 
\int_{\mathcal{R}_{\mathfrak{s}}(t)}
(\nabla \cdot  \zeta)
(\nabla \cdot  \eta_c) \ d x 
-
\int_{\Gamma_1(t)} 
 \sigma_f( v,p_f) \  n_c \cdot  \eta_c \ d\Gamma
\vspace{2mm} \\
\displaystyle
+
\int_{\Gamma_1(t)} 
 \sigma_c(\mathfrak{h}) \  n_c \cdot  \eta_c \ d\Gamma
=
\int_{\mathcal{R}_{\mathfrak{s}}(t)}
[\textbf{div} \  \sigma_c(\mathfrak{h})] \cdot  \eta_c \ d x
.
\end{cases}
\end{equation}
%%%%
%%%
%%%
Consider a time step $\Delta t>0$ and a finite element partition $\mathcal{U}_h$ of the solidification zone $\mathcal{R}_{\mathfrak{s}}(t)$ of maximum diameter $h$. Our aim is to approximate the solution $ \zeta$ at time $t_n=n \Delta t$, for $n \in \N$ in the finite element space. At any time $t$ consider the finite dimensional sub-space 
\begin{align*}
U_h
=
\{
 \eta_h:
 \eta_h
=
\eta_1 \psi_1
+
\ldots
+
\eta_N \psi_N
\}
\subset
\mathcal{W}_c,
\end{align*}
where $\{ \psi_i \}_i$ is a family of linearly independent functions with compact support, which are piecewise polynomials. In particular, we consider them to be of degree 2. Thus, at the instant $t=t_n$, the discretized formulation is
%%%%%%%
%%%%%%%
\begin{equation}\label{V.F-elasticity eqn}
\begin{cases}
\displaystyle
2 \mu_c
\int_{\mathcal{R}_{\mathfrak{s}}(t)}
 \varepsilon( \zeta_h^n):  \varepsilon( \eta_h) \ d x
+
\lambda_c 
\int_{\mathcal{R}_{\mathfrak{s}}(t)}
(\nabla \cdot  \zeta_h^n)
(\nabla \cdot  \eta_h) \ d x 
-
\int_{\Gamma_1(t)} 
 \sigma_f( v^n,p^n_f) \  n_c \cdot  \eta_h \ d\Gamma
\vspace{2mm} \\
\displaystyle
+
\int_{\Gamma_1(t)} 
 \sigma_c(\mathfrak{h}_h^n) \  n_c \cdot  \eta_h \ d\Gamma
=
\int_{\mathcal{R}_{\mathfrak{s}}(t)}
[\textbf{div} \  \sigma_c(\mathfrak{h}_h^n)] \cdot  \eta_h \ d x \qquad \forall \  \eta_h \in U_h
.
\end{cases}
\end{equation}
%%%%
%\begin{align}\label{V.F-elasticity eqn-discrete}
%2 \mu_s 
%\int_{\mathcal{R}_{\mathfrak{s}}(t)}
% \varepsilon( u^n_h):  \varepsilon( \eta_h) \ d x
%+
%\lambda_s 
%\int_{\mathcal{R}_{\mathfrak{s}}(t)}
%(\nabla \cdot  u^n_h)
%(\nabla \cdot  \eta_h) \ d x 
%-
%\int_{\Gamma_1(t)} 
% \sigma_f( v^n,p^n_f) \  n_c \cdot  \eta_h \ d\Gamma
%=0.
%\end{align}
%
%%%
%%%
Upon solving \eqref{V.F-elasticity eqn} using FreeFem++ software we obtain the displacement of the domain $\mathcal{R}_{\mathfrak{s}}(t)$, consequently, we get its deformation that is illustrated on Figure \ref{disp_coag}.
\\

\begin{figure}[h!]
  \centering
  \includegraphics[scale=0.4]{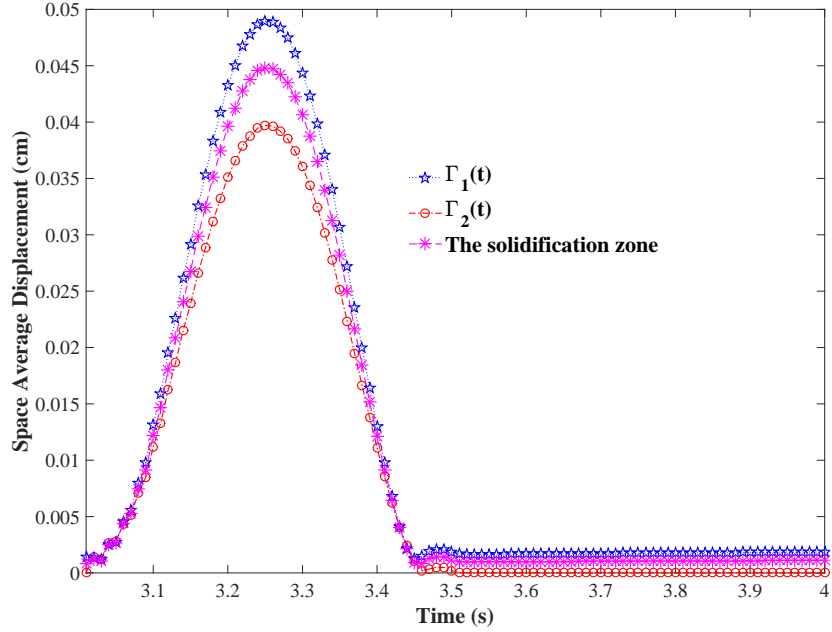}
\caption{The space average displacement of the solidification zone and its boundaries $\Gamma_1(t)$ and $\Gamma_2(t)$.}
\label{Displacement on boundaries}
\end{figure}

\begin{figure}[h!]
  \centering
  \includegraphics[scale=0.4]{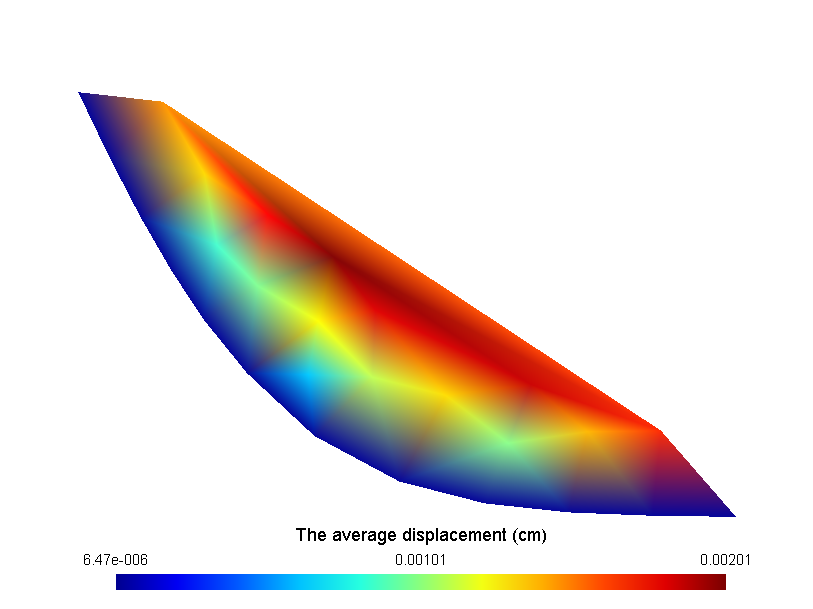}
\caption{The time average displacement of the solidification zone $\mathcal{R}_{\mathfrak{s}}(t)$.}
\label{av-disp}
\end{figure}

The deformation of the stenosis due to the blood external stress, results a deformation of the solidification zone. Indeed, this is due to the continuity of displacements on the stenosis-zone interface $\Gamma_2(t)$ given by the condition \eqref{Continuity of Displacement}.
The graphs corresponding to the space average displacement $\overline{ u}_{av}(t)$ of the solidification zone and the boundaries  $\Gamma_1(t)$ and $\Gamma_2(t)$ during the time interval between 3 s and 4 s, are illustrated on Figure \ref{Displacement on boundaries}. The space average displacement $\overline{ u}_{av}(t)$ is defined by
\begin{align*}
\overline{ u}_{av}(t)=
\dfrac{1}{|\mathcal{R}_{\mathfrak{s}}(t)|} 
\int_{\mathcal{R}_{\mathfrak{s}}(t)} \norm   u( x,t)   \norm_2 \ d  x ,
\end{align*} 
where $\norm . \norm_2$ is the Euclidean norm in $\R^2$ and $|\mathcal{R}_{\mathfrak{s}}(t)|$ is the area of $\mathcal{R}_{\mathfrak{s}}(t)$ provided that it is strictly positive. In a similar way we define the space average displacements of $\Gamma_1(t)$ and $\Gamma_2(t)$.
%%%
\\
From Figure \ref{Displacement on boundaries} we observe that the boundary $\Gamma_1(t)$ possesses the highest displacement. This fact is shown on Figure \ref{av-disp} which shows the time average displacement $\overline { u}_k( x)$ of the solidification zone on the time interval $3$ s-$4$ s.
The time average displacement at any position $ x$ is given by the formula
\begin{align*}
\overline{ u}_k( x)=
\dfrac{1}{k+1}
\sum_{i=0}^k  u( x,i),
\end{align*} 
where $ u( x,i)$ is the displacement of the solidification zone at any time $i \in \N^*$.
%%%
\\

It seems reasonable for $\Gamma_1(t)$ to possess the highest displacement, in fact, the displacement of the boundary $\Gamma_2(t)$ represents the displacement of the stenosis-zone boundary, while the displacement of the boundary $\Gamma_1(t)$ is a result of the deformation of the whole zone. The high displacement on $\Gamma_1(t)$ rises our curiosity to analyze the external stress exerted by blood on this boundary. Its average is given by the expression 
\begin{align*}
\dfrac{1}{|\Gamma_1(t)|} 
\int_{\Gamma_1(t)}  \sigma_f( v,p_f ) \  n_c \ d\Gamma,
\end{align*} 
where $|\Gamma_1(t)|$ stands for the length of the border $\Gamma_1(t)$ provided that its length is strictly positive. 
 \begin{figure}[h!]%{0.5\linewidth}
   \centering
   \includegraphics
   [
   %trim=2cm 1cm 4cm 1cm,height=5cm,
   scale=0.4]{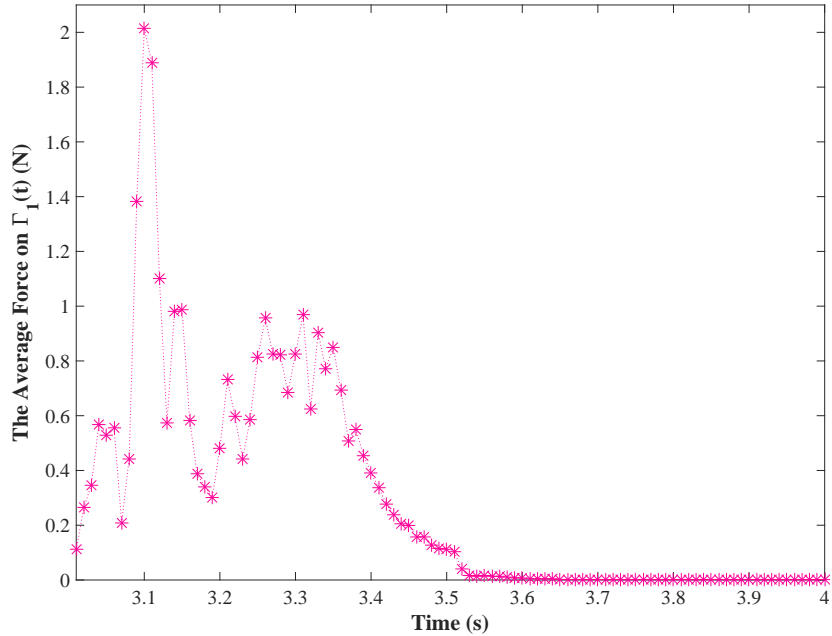}
    \caption{The magnitude of the average external force exerted on $\Gamma_1(t)$ at any time $t$.}
    \label{Force-magnitude}
  \end{figure}%%
\\

\begin{figure}[h!]%{0.5\linewidth}
   \centering
   \includegraphics
   [
   %trim=2cm 1cm 4cm 1cm,height=5cm,
   scale=0.4]{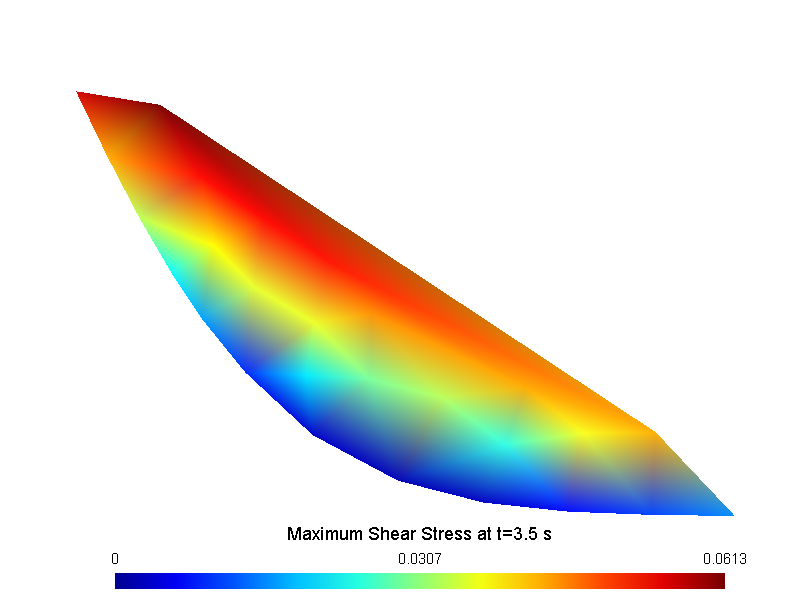}
    \caption{The maximum shear stress within the solidification zone at time $t=3.5$ s  (N/cm$^2$).}
    \label{Coag-Shst}
  \end{figure}

Figure \ref{Force-magnitude} shows that the magnitude of the average force exerted by the blood flow on the boundary $\Gamma_1(t)$ is large, hence, it results an inward resistance effect on this boundary which is large compared to the average displacement of the boundary $\Gamma_1(t)$ (see Figure \ref{Displacement on boundaries}). In other words, the force on the boundary $\Gamma_1(t)$ is opposed by the deformation of the solidification zone resulting from the deformation of the stenosis. Whence, the stress exerted on $\Gamma_1(t)$ will form a resistance factor against the displacements of $\Gamma_1(t)$ and $\Gamma_2(t)$, which will end up with the fragmentation of the crusted solidified blood.
To investigate the effect of the stress on the solidification zone, we will analyze the maximum shear stress $ \sigma_{max}$ given by the expression \eqref{max-shs}. Its pattern within the solidification zone 
is illustrated on Figure \ref{Coag-Shst} at $t=3.5$ s.
 %%
\\

\begin{figure}[h!]%{0.5\linewidth}
   \centering
   \includegraphics
   [
   %trim=2cm 1cm 4cm 1cm,height=5cm,
   scale=0.4]{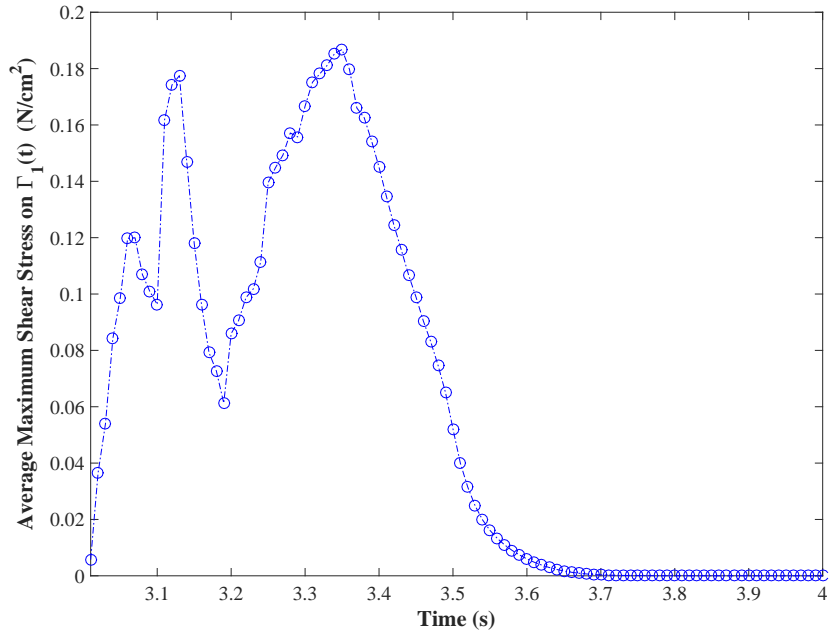}
    \caption{The magnitude of the average maximum shear stress on $\Gamma_1(t)$ at any time $t$.}
    \label{Maxshst-Gamma1}
  \end{figure}

Figure \ref{Coag-Shst} shows that the upper part of the solidification zone and the blood-zone interface $\Gamma_1(t)$ possesses the highest maximum shear stress value. Indeed, as we have mentioned previously, as we come closer to the peak of the stenosis its stiffness decreases, thus its displacement increases and it will deform easily, consequently, the upper part of the solidification zone will be more affected by the displacement at the interface $\Gamma_2(t)$. Further, the stress exerted by blood on $\Gamma_1(t)$ which is of high magnitude (see Figure \ref{Force-magnitude}) will lead to a high maximum shear stress.
The space average maximum shear stress on the boundary $\Gamma_1(t)$ at any time $t$ is given by the formula 
\begin{align*}
\overline{\sigma}_{max}(t)=
\dfrac{1}{|\Gamma_1(t)|} 
\int_{\Gamma_1(t)} \norm  \sigma_{max} ( x,t)  \norm_2 \ d  x ,
\end{align*} 
The graph corresponding to the magnitude of the average maximum shear stress $ \overline{\sigma}_{max}$ on the boundary $\Gamma_1(t)$ is given in Figure \ref{Maxshst-Gamma1}.
  \\
   \\
   
The force exerted by the blood on the solidification zone which opposes its displacement will form a frictional force on the solidification zone (digging manner). Further, from Figure \ref{Maxshst-Gamma1}, we observe that at the instant $t=3.5$ s when the solidification zone returns to its equilibrium position, the boundary $\Gamma_1(t)$ is still under the impact of the maximum shear stress. As a result, the maximum shear stress will scrape the crust leading to the release of some pieces into the flow, which will block the flow at some levels of narrow vessels or arterioles causing an infarction.
\newpage
\section{Conclusion}
This work is devoted for proposing a mathematical model for the rupture of blood in stenosed arteries. A fluid-structure interaction problem representing the interaction between blood flow and an existing stenosis in arteries is considered. The blood is assumed to be a homogeneous non-Newtonian incompressible fluid whose dynamics is given by the incompressible Navier-Stokes equations and of a viscosity $\mu$ obeying Carreau model, while the arterial wall is a non-linear hyperelastic material described by the quasi-static elasticity equations.
The simulations have shown a deep view of what is happening in the stenosed artery and how it would be affected with some variables that we can analyze. They helped us in configuring the existence of mainly three remarkable regions (see Figure \ref{Three zones}).
\\

%Based on these results, we assume that the blood is a non-Newtonian fluid with a time-dependent viscosity $\mu$ that obeys Carreau model. 
We believe that what is ruptured is not the stenosis plaque, rather, the solidified blood near the stenosis. In fact, the fibrous cap is a stiffened part of the artery wall which is enlarged due to the inflammation beneath it, which rebut the assumption of being released into the flow. Hence, a first  step towards a rupture model is to locate the solidification zone as we believe that the jelly-like material in this zone is the ruptured substance. For this sake, the viscosity $\mu$ is reformulated so that it is a time-dependent function related to its history represented by the local in time average viscosity $\hat{\mu}_k$. Indeed, a transit from a liquid state to a jelly-like material is linked to an increase in the viscosity. 
\\

Based on the properties of viscous materials, we can assume that a solidification zone is characterized by a high viscosity and a negligible speed. By investigating the pattern of the average speed and the average viscosity, we consider a viscosity threshold $\mu_{th}$ such that when the computed blood viscosity $\mu$ exceeds it, regions $\mathcal{D}_\mu$ of high viscosity are detected. Similarly, if the speed of blood is less than the speed threshold $v_{th}$, then we locate the regions $\mathcal{D}_v$ possessing negligible speed. The solidification zone $\mathcal{R}_{\mathfrak{s}}$, spotted at the edge of the stenosis, is the intersection of the regions $\mathcal{D}_\mu$ and $\mathcal{D}_v$.
\\

A rupture model is derived based on the forces acting on the solidification zone $\mathcal{R}_{\mathfrak{s}}$. For this sake, semi-solidified blood is considered to be a linear elastic material that obeys Hooke's law and that is under the effect of an external stress from the blood and the deformation of the stenosis. Upon solving numerically the elasticity equations, results have showed an inward force resulting from the shear stress exerted by the blood on this zone, opposed by the zone deformation  due to the deformation of the stenosis. These opposite effects will lead to the fragmentation of the solidified blood of the solidification zone. Further, the maximum shear stress will scrape the crust of this zone. As a consequence, detached pieces will be drifted by the flow and at some sites will block the artery causing an infarction.
\\ 
\\

%%%%%%%%%%%%%%%%%%%%%%%%%%%%%%%%%%%%%%%%%%%%%%%%%%%%%%%%%
%%%%%%%%%%%%%%%%%%%%%%%%%%%%%%%%%%%%%%%%%%%%%%%%%%%%%%%%%
%%%%%%%%%%%%%%%%%%%%%%%%%%%%%%%%%%%%%%%%%%%%%%%%%%%%%%%%%
%%%%%%%%%%%%%%%%%%%%%%%%%%%%%%%%%%%%%%%%%%%%%%%%%%%%%%%%%

\noindent {\bf Acknowledgments.} The authors would like to thank the Rectorat of the Lebanese University for funding the research project about mathematical analysis and numerical simulation of the blood flow in stenosed arteries.

\end{document}